# $k$-Submodular Interdiction Problems under Distributional Risk-Receptiveness and Robustness: Application to Machine Learning

Seonghun Park[1] and Manish Bansal[2*]

[1,2]Grado Department of Industrial and Systems Engineering, Virginia Tech, Blacksburg, 24061, VA, USA.

*Corresponding author(s). E-mail(s): bansal@vt.edu;
Contributing authors: seonghun@vt.edu;

**Abstract**

We study submodular optimization in adversarial context, applicable to machine learning problems such as feature selection using data susceptible to uncertainties and attacks. We focus on Stackelberg games between an attacker (or interdictor) and a defender where the attacker aims to minimize the defender's objective of maximizing a $\boldsymbol{k}$-submodular function. We allow uncertainties arising from the success of attacks and inherent data noise, and address challenges due to incomplete knowledge of the probability distribution of random parameters. Specifically, we introduce Distributionally Robust $\boldsymbol{k}$-Submodular Interdiction Problem (DRO $\boldsymbol{k}$-SIP) and Distributionally Risk-Receptive $\boldsymbol{k}$-Submodular Interdiction Problem (DRR $\boldsymbol{k}$-SIP) along with finitely convergent exact algorithms for solving them. When solving the DRO $\boldsymbol{k}$-SIP, the attacker optimizes their expected payoff with respect to the worst-case probability distribution within the ambiguity set, and thereby have robust attack strategies despite distributional ambiguity. In contrast, the DRR $\boldsymbol{k}$-SIP identifies attacker strategies with the best-case probability distribution, and identifies critical vulnerabilities for the defender. The optimal values derived from both DRO $\boldsymbol{k}$-SIP and DRR $\boldsymbol{k}$-SIP offer a confidence interval-like range for the expected value of the defender's objective function, capturing distributional ambiguity. We conduct computational experiments on instances of feature selection and sensor placement problems, using Wisconsin breast cancer data and synthetic data, respectively.

# 1 Introduction

Submodularity is an important concept in the domain of combinatorial optimization as submodular functions encompass various classes of functions, including weighted coverage, entropy, and mutual information functions [1]. Consequently, the submodular functions have found extensive use in machine learning with applications such as feature selection [2, 3], image segmentation [4], data summarization [5], influence maximization [6], and sensor placement [7–9]. Recently, it has been observed that machine learning models are susceptible to unexpected noise associated to the training or testing data, and to the adversarial attacks where an adversary can inject deteriorated training data, thereby leading to inaccuracies in the output results [10]. For an example, the feature selection process becomes significantly more complex in adversarial scenarios due to the actions of attackers aiming to degrade model performance through feature removal or label manipulation [11].

In this paper, we consider min-max interdiction problem with $k$-submodular functions that involve two non-cooperating players: an interdictor (or attacker) and a follower (or defender). The attacker with limited budget aims to minimize the defender's objective of maximizing a submodular function by interdicting a subset of the ground set. The attacker's optimal solution identifies the subset of ground set whose disruption would most significantly undermine defender's possible objective and this approach enables the attacker to prioritize the resources effectively, by targeting the most vulnerable components for maximum disruption. We explore this framework across two distinct problems: the Feature Selection Interdiction Problem (FSIP) and the Weighted Coverage Interdiction Problem (WCIP), both characterized by the defender's objective function being submodular. In FSIP, the attacker targets specific data features to degrade the model performance and, as a response, a defender selects an optimal subset from the remaining features that maximizes performance of machine learning predictive models. Similarly, in WCIP, the attacker interdicts certain locations to prevent sensor installation, and the defender then places sensors at the remaining sites to maximize coverage. The objective function of the defender's problem in FSIP and WCIP are defined using submodular functions; refer to Sections 3.2 for details.

In this paper, we introduce stochastic submodular interdiction problems to handle uncertainties in the success of attacks and data associated with the ground set (e.g., uncertain data corresponding to each feature in a feature selection problem). Specifically, the attacker/interdictor aims to remove components (i.e., features or potential sensor locations) of the ground set to minimize the expected value of the defender's objective function. However, in many applications, the evaluation of expectation is constrained by the limited availability of historical data, which complicates the estimation of true probability distributions associated with uncertain data parameters. To address this additional challenge, we employ distributionally robust optimization (DRO) framework that optimizes objective function for the worst-case probability distribution within a predefined set, known as the ambiguity set [12–17]. This framework yields a distributionally robust solution for a risk-averse attacker, thereby offering a robust attacking strategy when the attacker is the main protagonist. From the defender's perspective, the assumption that the attacker's strategy is solely risk-averse could lead to a misjudgment of the attacker's desire for aggressive actions (receptiveness to take risk), potentially leaving the defender unprepared for more risky and impacting tactics. To this end, we introduce distributionally robust submodular interdiction problem (DRO SIP) and a distributionally risk-receptive submodular interdiction problem (DRR SIP). An optimal solution for the DRR SIP emphasizes key vulnerabilities, identifying which elements, if compromised, would most severely degrade their objectives, thereby guiding the development of targeted defensive strategies.

Furthermore, in many practical applications, the challenge often involves selecting multiple disjoint subsets instead of a single subset. For instance, consider the problem where defender aims to maximize area coverage by deploying $k \in \mathbb{Z}_+$ distinct types of sensors, each with a unique range, and are constrained to install at most one sensor for each location. This situation underscores the need for a mathematical model that can address the selection of multiple, distinct subsets to optimize a given objective. To incorporate this feature, we employ $k$-submodular function [18], which reduces to standard submodular function when $k = 1$, as a defender's objective function in stochastic SIP, DRO SIP and DRR SIP. We denote these generalized problems as stochastic $k$-SIP, DRO $k$-SIP and DRR $k$-SIP, and present exact solution approaches to solve them. This incorporation not only addresses the problem of placing $k$-types of sensors for maximizing the coverage but also extends applicability of our research to any situation where $k$-submodular function can be employed [19, 20].

## 1.1 Contributions and Organization of the Paper

As per our knowledge, stochastic $k$-SIP, DRO $k$-SIP and DRR $k$-SIP have not been addressed in current literature for any $k \geq 1$ and moreover, deterministic $k$-SIP has not been studied for $k \geq 2$. We review the literature related to these problems in Section 2. In Section 3, we present necessary background for submodular and $k$-submodular functions and introduce the formulations of deterministic $k$-SIP, stochastic $k$-SIP, DRO $k$-SIP and DRR $k$-SIP. Next, in Section 4, we describe exact solution methodologies for solving DRO $k$-SIP and DRR $k$-SIP (that subsume both deterministic and stochastic $k$-SIPs), by introducing families of valid inequalities and embedding them within decomposition-based approaches. In Sections 5 and 6, we present results of our computational experiments and concluding remarks, respectively.

To evaluate the impact of this paper, we use the Wisconsin Breast Cancer Data [21] for FSIP and obtain optimal solutions for both attacker (interdicted features) and defender (selected features) by solving deterministic, DRO, stochastic, and DRR $k$-SIPs for $k = 1$. Four Support Vector Classifiers (SVCs) were trained using only the defender-selected features from each problem types and tested against 100 noise-injected scenarios, to estimate the robustness of attacking strategies. Figure 1 illustrates the performance of SVCs, where each dot represents the test accuracy of the SVCs for noise-injected data sets, and the enveloping colored curves represent the frequency of occurrence of these accuracy dots out of the 100 scenarios and two flat lines in the curvature represent the mean (along with mean values) and the median (without mean values). Observe that the test accuracy range for the SVC trained with DRO-SIP solution aligns closely with those trained using solutions from both deterministic and stochastic (risk-neutral) SIPs, i.e., from 89% to 95%. However, the test accuracy range for the SVC trained with the solution from DRR SIP spans from 82% to 96%, presenting a broader distribution compared to other attacking strategies. This suggests that attacking strategy derived from DRR-SIP has the potential to degrade prediction model more significantly compared to other attacking strategies. In a nutshell, this analysis provides not only bounds on the accuracy for varying levels of risk-appetite (ranging from risk-aversion to risk-receptiveness) of decision makers but also returns most vulnerable features of the data set. (Refer to Section 5.2 for results of computational experiments conducted on randomly generated instances of weighted coverage problem.)

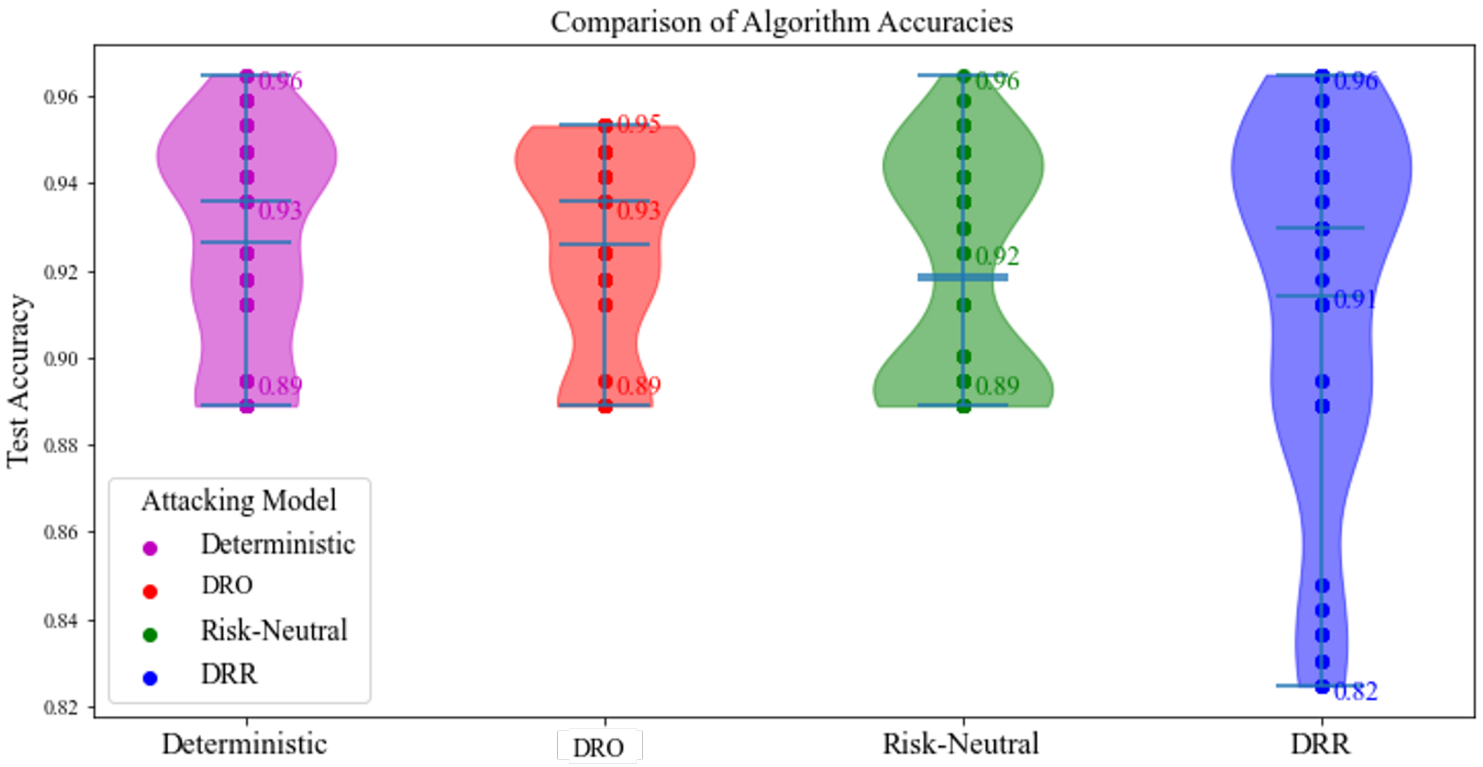

**Fig. 1**: Performance of Support Vector Classifiers (SVCs) on Wisconsin Breast Cancer Data [21] using defender-selected features provided by different modeling framework

## 2 Literature Review

In this section, we review the literature related to submodular function, $k$-submodular function, Stackelberg zero-sum games (also referred as interdiction games), and adversarial machine learning including feature selection problems.

### 2.1 $k$-Submodular Functions and Their Applications

Submodular functions have a wide range of applications. Krause and Golovin [1] provided a comprehensive survey of the applications of submodularity in machine learning. Vohra and Hall [22] solved maximal covering problem by leveraging submodular objective function. Likewise, Nemhauser et al. [23] showed that the objective function of the incapacitated facility location problem can be formulated as a monotone submodular function and provided a greedy algorithm with $(1-1/e)$ approximation ratio for solving it. In another direction, Wei et al. [3] introduced framework for unsupervised data subset selection problem by formulating maximum likelihood estimator function of Naïve Bayes classifier and Nearest Neighbor classifiers as submodular function and selecting the data subset which maximizes this submodular function using the greedy approach. Recently, Kothawade et al. [24] introduced various submodular information measures for data subset selection with the aim of selecting subsets with desired characteristics.

A generalization of the submodular function, referred to as a $k$-submodular function, was introduced by Huber and Kolmogorov [18]. As mentioned earlier, it has been widely applied to the applications such as multi-topic influence maximization problem [19, 25], multi-type sensor placement problem [19, 26] and information coverage problem [26]. Acknowledging its practicability, most of the research focus on the methodology to maximize $k$-submodular functions. Ohsaka and Yoshida [19] introduced a greedy algorithm for monotone $k$-submodular function maximization under two types of size constraints: an overall cardinality constraint and individual cardinality constraints for each of the $k$ subsets. Their greedy approach achieved approximation ratios of $1/2$ and $1/3$ for these respective constraints. Expanding upon this, Sakaue [27] addressed matroid constraints, yielding a $1/2$-approximation algorithm, while Tang et al. [28] focused on non-negative monotone $k$-submodular function maximization under knapsack constraints, attaining an approximation ratio of $\left(\frac{1}{2} - \frac{1}{2e}\right)$. While the majority of research has been focused on approximation methods, Yu and Küçükyavuz [29] offered an exact solution approach. They first introduced valid inequality which is the tight upper approximation of the hypograph of any $k$-submodular function and incorporate them in the delayed constraint generation method for repeatedly refining hypograph of $k$-submodular function to solve the problem exactly.

## 2.2 Stackelberg Zero-Sum Games

The attacker-defender dynamics is commonly modeled through Stackelberg zero-sum games. A representative problem of this domain is the Network Interdiction Problem (NIP), [30, 31], where an attacker aims to disrupt a network of nodes and arcs, thereby hindering a defender's (also known as an evader or follower) ability to transport illegal drugs or nuclear materials. The goal of the defender is to either maximize the flow [32] or minimize the shortest path [31] from source node to destination node, and the attacker's aim is to minimize or maximize, respectively, the defender's objective. Israeli and Wood [33] introduced stochastic network interdiction problem by incorporating uncertain data parameters defined by random variables with known probability distribution. Readers can refer to [34] for a comprehensive survey on other variants of NIP. Building on the literature, Kang and Bansal [35] considered incomplete information of probability distribution associated with uncertain parameters in stochastic NIP, and offered adjustments based on varying level of risk-appetite of decision makers in NIP. In another direction, Park and Bansal [36] extended beyond the conventional formulations of the defender's problem as linear or integer programs. They considered a computational geometry problem of placing multiple rectangular camera view-frames to capture areas with maximum threat as defender's problem. Similarly, Tanınmış and Sinnl [37] explored a deterministic interdiction game with a monotone 1-submodular function as the defender's objective and presented a branch-and-cut based exact solution approach for solving this problem.

### 2.3 Machine Learning in the Presence of Adversarial Attacks

For the sake of completeness of literature review, we briefly review some advances in adversarial machine learning that are relevant for this paper. Kurakin et al. [38] pointed that machine learning models are susceptible to the adversarial attacks such as intentionally injected malicious input samples or noise. They emphasized the importance of constructing robust models that perform well in various adversarial scenarios. For feature selection problem, Globerson and Roweis [39] introduced a problem where a subset of features of test instances is interdicted and thereby the support vector machine model trained using labeled train data with all features, can only rely on the subset of features associated with test data for the classification. This led to a game-theoretic min-max problem allowing defender (i.e., data user) to build a model that is resilient to feature interdiction executed by the attacker. Dekel and Shamir [40] presented a more efficient algorithm for solving the foregoing problem by taking into account the level of importance of each feature. Both studies optimized the objective function based on the labeled training data, thereby relying provided label information. They developed robust models for test samples, where only a subset of the training data features is available, and different features may be interdicted for different test data points. However, they did not account for uncertainty and incomplete information regarding the probability distribution. Additionally, their approaches depend on the label information, which can be vulnerable to attacks.

In the domain of Graph Neural Networks, addressing adversarial attack on graph structures has become a critical area of research. Dai et al. [41] proposed reinforcement learning-based effective attacking strategy by modifying or deleting edges of the graph structure where the model learns to attack the graph from the classifier's prediction. Xu et al. [42] also introduced a gradient-based attack method to obstruct the processing of graph structured data by data users, and adversarial learning techniques for graph neural network to build more robust models against such attacks. As per our knowledge, the aforementioned studies do not incorporate the following features: uncertainty in the data, incomplete information of probability distribution, and adjustments based on risk-appetite of a decision maker.

## 3 $k$-Submodular Interdiction Problems without and with Uncertainty: Formulations and Applications

In this section, we provide some definitions necessary to introduce formulations of four problems: Deterministic, Stochastic (Risk-Neutral), Distributionally Robust (DRO), and Risk-Receptive (DRR) $k$-SIPs. We also present two applications of these formulations, i.e., feature selection interdiction problem and weighted coverage interdiction problem, that we employ for our computational experiments as well.

**Definition 1** (Submodular Function)**.** *Let $N = \{1, \dots, n\}$ be a non-empty finite ground set, and let $2^N$ denote the power set of $N$, i.e., the set of all subsets of $N$. A function $f : 2^N \to \mathbb{R}$ is a monotone submodular if for every $A \subseteq B \subseteq N$ and $i \in N \setminus B$,* *the following diminishing returns property holds:*

$$\rho_i(A) = f(A \cup \{i\}) - f(A) \geq f(B \cup \{i\}) - f(B) = \rho_i(B), \tag{1}$$

*where $\rho_i(A)$ and $\rho_i(B)$ denotes the marginal gain of adding element $i$ to set $A$ and $B$, respectively.*

**Definition 2** (Set of $k$ Disjoint Sets)**.** *For $k \geq 1$ and $N = \{1, \ldots, n\}$, we denote a set of all $k$ disjoint subsets of $N$ by*

$$\mathbb{X}(N, k) = \big\{(Z_1, \ldots, Z_k) : Z_q \subseteq N \textit{ for all } q \in \{1, \ldots, k\} \textit{ and } Z_q \cap Z_{q'} = \emptyset \text{ for } q \neq q'\big\}.$$

**Definition 3** ($k$-Submodular Function)**.** *A function $f : \mathbb{X}(N, k) \to \mathbb{R}$ is a $k$-submodular function if for any $\mathbf{X} = (X_1, \ldots, X_k)$ and $\mathbf{Y} = (Y_1, \ldots, Y_k) \in \mathbb{X}(N, k)$, the following holds:*

$$f(\mathbf{X}) + f(\mathbf{Y}) \geq f(\mathbf{X} \sqcap \mathbf{Y}) + f(\mathbf{X} \sqcup \mathbf{Y}) \tag{2}$$

*where $\mathbf{X} \sqcap \mathbf{Y} = (X_1 \cap Y_1, \ldots, X_k \cap Y_k)$ and*

$$\mathbf{X} \sqcup \mathbf{Y} = \left\{ (X_i \cup Y_i) \setminus \bigcup_{\substack{q=1 \\ q \neq i}}^{k} (X_q \cup Y_q), \right\}_{i=1}^{k}$$

**Definition 4** (Marginal Gains for $k$-Submodular Functions)**.** Given $\mathbf{X} \in \mathbb{X}(N, k)$, for any $q \in \{1, \ldots, k\}$, $i \in N \setminus \bigcup_{r=1}^{k} X_r$, the marginal gain of adding $i$ to $X_q$ is defined as

$$\rho_{q,i}(\mathbf{X}) = f(X_1, \ldots, X_q \cup \{i\}, \ldots, X_k) - f(\mathbf{X}).$$

If function $f$ is $k$-submodular, then for every $\mathbf{X}, \mathbf{Y} \in \mathbb{X}(N, k)$ such that $X_r \subseteq Y_r$ for all $r \in \{1, \ldots, k\}$,

$$\rho_{q,i}(\mathbf{X}) \geq \rho_{q,i}(\mathbf{Y}) \qquad \textit{for all } i \in N \setminus \bigcup_{r=1}^{k} Y_r \textit{ and } q \in \{1, \ldots, k\}.$$

**Definition 5** ($k$-Submodular Maximization)**.** *Given a set of problem dependent feasible solutions $\mathcal{S} \subseteq \mathbb{X}(N, k)$, the $k$-submodular maximization problem is defined as:*

$$\max\{f(\mathbf{S}) : \mathbf{S} \in \mathcal{S}\}. \tag{3}$$

*For any $\mathbf{S} \in \mathbb{X}(N, k)$, we can write $f(\mathbf{S}) = f(\mathbf{s})$ where $\mathbf{s} = (s_1, \ldots, s_k) \in \{0, 1\}^{kn}$ and $s_{q,i} = 1$ when $i \in S_q$, and $s_{q,i} = 0$ otherwise for all $q \in \{1, \ldots, k\}$ and $i \in N$.*

## 3.1 Formulations for *k*-SIPs without and with Uncertainty

To illustrate the formulation of the deterministic $k$-submodular interdiction problem, we consider an example of the weighted coverage interdiction problem. In this problem, an interdictor intends to prevent a defender from installing sensor of type $q \in \{1, \ldots, k\}$ at location $i \in N$ by minimizing the defender's objective of maximizing coverage. The interdictor's decision variable $x_{q,i} = 1$ if he decides to block the installation of a type $q$ sensor at sites $i$ and 0 otherwise. In contrast, we denote the defender's decisions by $\mathbf{S} = (S_1, \ldots, S_k) \in \mathbb{X}(N, k)$ where $S_q$ represents a set of locations where the defender choose to place type $q \in \{1, \ldots, k\}$ sensors. We assume that both players are subject to budget limitations for each sensor type $q \in \{1, \ldots, k\}$, denoted as $A_q$ and $D_q$ for the interdictor and defender, respectively. We formulate the deterministic $k$-submodular interdiction problem as:

**Deterministic $k$-SIP**

$$\min_{\mathbf{x} \in \mathcal{X}} \quad \Phi_D(\mathbf{x}) \tag{4a}$$

$$\text{where } \mathcal{X} := \left\{ \mathbf{x} \in \{0,1\}^{kn} : \sum_{i=1}^{n} x_{q,i} \leq A_q, \ \forall q \in \{1, \ldots, k\}, \sum_{q=1}^{k} x_{q,i} \leq B_i, \ \forall i \in N \right\}, \tag{4b}$$

$$\Phi_D(\mathbf{x}) := \max \left\{ f(\mathbf{S}) : \mathbf{S} \in \mathbb{X}(N, k), s_{q,i} \leq 1 - x_{q,i}, \ \forall q \in \{1, \ldots, k\}, \forall i \in N, \right. \tag{4c}$$

$$\left. \sum_{i=1}^{n} s_{q,i} \leq D_q, \ \forall q \in \{1, \ldots, k\} \right\}. \tag{4d}$$

In (4b), first constraint enforces the budget limits for the interdictor for each $q \in \{1, \ldots, k\}$, and the second constraint ensures that the interdictor selects $k$ disjoint subsets, i.e., for each location $i \in N$, an interdictor can block the installation of at most $B_i$ sensors of different types. The constraint in (4c) indicates that the defender cannot place sensor type $q \in \{1, \ldots, k\}$ at location $i \in N$ if it has been blocked by the interdictor's solution $\mathbf{x}$. Constraint (4d) enforces defender's budget limits for each $q \in \{1, \ldots, k\}$.

To incorporate uncertainty into problem (4), we introduce two random variables: a parameter $\mathcal{D}$ that affects the defender's $k$-submodular objective function, and a binary vector $\xi \in \{0,1\}^n$ representing the success of attacks on elements of the ground set $N$. Let $\Omega := \{\omega_1, \ldots, \omega_{|\Omega|}\}$ be a finite set of possible realizations of $(\xi, \mathcal{D})$, and each scenario $\omega \in \Omega$ corresponds to a realization $(\xi^\omega, \mathcal{D}^\omega)$ with associated probability $p_\omega$.

The realization $\mathcal{D}^\omega$ governs the marginal contributions of elements in the objective, while $\xi^\omega$ determines which elements are successfully interdicted and thus modifies the defender's feasible set. Specifically, we denote the scenario-dependent submodular function as $f^\omega(\cdot) := f(\cdot; \mathcal{D}^\omega, \xi^\omega)$, where the functional form is fixed (as a $k$-submodular function in our case), but its evaluation varies across scenarios due to the realization of the random parameters. Moreover, contrary to the assumption in problem (4) that all attacks are successful, we now assume that an attack on location $i \in N$ for sensor type $q \in \{1, \ldots, k\}$ is successful in scenario $\omega \in \Omega$ if and only if $\xi_i^\omega \cdot x_{q,i} = 1$, where $\xi^\omega = (\xi_1^\omega, \ldots, \xi_n^\omega)$. After observing the attacker's solution $\mathbf{x} \in \mathcal{X}$ and a realization of uncertainties $\omega \in \Omega$, the defender selects subset of available items, denoted by $\mathbf{S}^\omega$, to maximize its $k$-submodular function. For $\omega \in \Omega$, the defender's problem is defined as

$$Q_\omega(\mathbf{x}, \xi^\omega) := \max \Bigg\{ f^\omega(\mathbf{S}^\omega) : \mathbf{S}^\omega \in \mathbb{X}(N, k), \tag{5a}$$

$$s^\omega_{q,i} \leq 1 - x_{q,i}\xi^\omega_i, \ \forall q \in \{1, \dots, k\} \text{ and } i \in N, \tag{5b}$$

$$\sum_{i=1}^{n} s^\omega_{q,i} \leq D_q, \ \forall q \in \{1, \dots, k\} \Bigg\}. \tag{5c}$$

If probability distribution $P$ associated with random parameters is known, we get **Stochastic (Risk-Neutral) $k$-SIP**

$$\min_{\mathbf{x} \in \mathcal{X}} \ \Phi_N(\mathbf{x}) := \mathbb{E}_P[Q_\omega(\mathbf{x}, \xi^\omega)], \tag{6}$$

where $\mathbb{E}_P[\cdot]$ is an expectation operator.

When complete information of the probability distribution $P$ is not known and it belongs to a predefined ambiguity set $\mathfrak{P}$, we formulate Distributionally Risk-Receptive $k$-Submodular Interdiction Problem (DRR $k$-SIP) and Distributionally Robust $k$-Submodular Interdiction Problem (DRO $k$-SIP) with the predefined ambiguity set $\mathfrak{P}$ as follows.

**DRR $k$-SIP**

$$\min_{\mathbf{x} \in \mathcal{X}} \ \left\{ \Phi_R(\mathbf{x}) := \min_{P \in \mathfrak{P}} \mathbb{E}_P[Q_\omega(\mathbf{x}, \xi^\omega)] \right\} \tag{7}$$

**DRO $k$-SIP**

$$\min_{\mathbf{x} \in \mathcal{X}} \ \left\{ \Phi_A(\mathbf{x}) := \max_{P \in \mathfrak{P}} \mathbb{E}_P[Q_\omega(\mathbf{x}, \xi^\omega)] \right\} \tag{8}$$

where the inner minimization and maximization with respect to the ambiguity set handle the risk-receptiveness and risk-aversion of the attacker by selecting best-case and worst-case distributions, respectively. In case $|\mathfrak{P}| = 1$, formulations (7) and (8) reduce to (6). Also, when $|\Omega| = 1$, formulations (6), (7), and (8) reduce to (4).

**Remark 1.** The general formulations for distributionally risk-receptive $k$-SIP (7) and distributionally robust $k$-SIP (8) can accommodate any choice of ambiguity set $\mathfrak{P}$ such as Wasserstein metric [43], moment information [44], $\phi$-divergence [45], and many more, and the solution approaches presented in this paper are applicable for these ambiguity sets. In Section 4, we provide sufficient conditions for the ambiguity sets under which the proposed approaches are finitely convergent.

## 3.2 Applications of $k$-SIPs

### 3.2.1 Feature Selection Interdiction Problem: An Application of 1-SIP

Feature selection is a fundamental problem in machine learning, where the goal is to identify a subset of relevant features that contribute most significantly to the predictive performance of a machine learning model. This process not only enhances model interpretability but also improves efficiency by reducing dimensionality and mitigating overfitting issues. Within this context, this paper addresses a specific challenge associated to the feature selection problem in the presence of an adversary. Specifically, data sets are often compromised through targeted attacks on specific features, thereby impacting the performance of predictive models employed by defenders. In our formulation, the attacker strategically selects a subset of features to be interdicted (e.g., corrupted or blocked) before the defender acts. This setting is modeled as an interdiction game, where the attacker acts as the leader and anticipates how the defender will respond. After observing the attack, the defender selects an unaffected subset of features to train the model, aiming to mitigate the degradation in predictive performance. Submodular functions are employed to assist in this feature selection process, a method proven effective by [2, 46] for choosing subset of features. This strategy is centered around the use of selected features for model training and focusing on maintaining model efficacy rather than enhancing performance directly.

Mathematically, FSIP is defined as follows. Let $N = \{1, \ldots, n\}$ be the index set of features, and let $X_i$ and $X_j$ denote the $i$-th and $j$-th feature vectors in the dataset $\mathcal{D} \in \mathbb{R}^{m \times n}$, each represented as a length-$m$ vector over $m$ data points. The similarity score $w_{i,j} (\geq 0)$ between features $i$ and $j$ is computed using a measure such as mutual information or cosine similarity between $X_i$ and $X_j$, depending on the application and data type. Then, function $f : 2^N \to \mathbb{R}$ defined as

$$f(S) = \sum_{i \in N} \max_{j \in S} w_{i,j}$$

is a monotone submodular function [2]. This function evaluates a given set $S \in N$ by selecting, for each $i \in N$, a feature $j \in S$ that has the highest similarity to $i$ as given by $w_{i,j}$, and then aggregating these values. Using this function, we solve problems (6), (7), and (8) and obtain optimal solutions for an attacker and defender. For details on the experimental setup and computational results, please refer to section 5.1.

### 3.2.2 Weighted Sensor Coverage Interdiction Problem as $k$-SIP

Consider a set $N = \{1, \ldots, n\}$ of potential sites for installation of sensors that are of $k$ different types depending on their monitoring ranges. Let $\mu_{i,q}$ be the reward obtained by covering site $i \in N$ using sensor type $q \in \{1, \ldots, k\}$. However, only one sensor can be placed at each site, and at most $D_q$ sensors of type $q \in \{1, \ldots, k\}$ can be installed. Let $\mathbf{S} = (S_1, S_2, \ldots, S_k) \in \mathbb{X}(N, k)$ denotes a feasible sensor placement strategy with $S_q$ representing the set of locations where type $q$ sensors are placed. Define $C_i = \{q \mid i \in S_q\}$ as the set of sensor types covering site $i$. Then the total reward is given by the $k$-submodular function $f : \mathbb{X}(N, k) \to \mathbb{R}$, i.e.,

$$f(\mathbf{S}) = \sum_{i \in N} \max_{q \in \mathbf{S}(i)} \mu_{i,q}.$$

In the interdiction framework, the attacker first selects a blocking strategy $\mathbf{x} \in \mathcal{X}$ to disable specific sensor placements. Given this interdiction, the defender chooses a placement strategy $\mathbf{S} = (S_1, \ldots, S_k)$ that maximizes the $k$-submodular reward $f(\mathbf{S})$ under uncertainty. The attacker's goal is to reduce the defender's utility by anticipating this optimal adaptation. Again, using this function, formulations (6), (7), and (8) lead to weighted sensor coverage interdiction problem under uncertainty. In Section 5.2, we describe the experimental setup and results in details.

# 4 Solution Methodologies: Valid Inequalities and Algorithms

We present valid inequalities for deterministic, DRO $k$-SIP and DRR $k$-SIP in Section 4.1, and decomposition-based exact approaches for solving DRO $k$-SIP and DRR $k$-SIP in Section 4.2. Specifically, we introduce tight lower approximations for the objective functions of problems (4), (6), (7) and (8) using affine functions that lead to inequalities, which are referred to as valid inequalities (or cutting planes). These inequalities are used to iteratively refine the lower bounds of the objective functions within a decomposition framework, thereby ensuring a global optimal solution (4.2). Note that these objective functions $\Phi_D(\mathbf{x})$, $\Phi_N(\mathbf{x})$, $\Phi_R(\mathbf{x})$, and $\Phi_A(\mathbf{x})$ are non-convex, non-concave, and non-submodular. This complicates the application of known optimization techniques, such as gradient methods or submodular cuts, to optimize them.

Recall that in this paper we assume a finite number of scenarios, i.e. $|\Omega| < \infty$. In the literature, problems with continuous support set have also been studied, where the second-stage problem is defined as a linear or mixed-integer linear program [47] and duality-based approaches are used to derive convex approximations or reformulations [43]. However, these approaches are dependent on the selection of the ambiguity set and lead to different reformulations. In contrast, the approaches presented in this paper are applicable to any general ambiguity set.

## 4.1 Valid Inequalities

The term "valid" in the valid inequality implies that all feasible solutions satisfy this inequality, and hence, the inequality is valid. In the following theorems, we present these inequalities and, for the sake of ease for the reader, all validity proofs are provided in Section 4.3. These inequalities provide a lower bounding approximation of the non-submodular function using affine function that we use within our decomposition algorithms for solving deterministic, stochastic, DRO, and DRR $k$-SIPs.

**Theorem 1.** *Given an attacker's solution $\hat{\mathbf{x}}$ and associated defender's optimal solution $\hat{\mathbf{S}}$, i.e., $\Phi_D(\hat{\mathbf{x}}) = f(\hat{\mathbf{S}})$ in Formulation (4), the following inequality is valid for all*

$\mathbf{x} \in \mathcal{X}$:

$$\Phi_D(\mathbf{x}) \geq \Phi_D(\hat{\mathbf{x}}) - \sum_{q=1}^{k} \sum_{i \in \hat{S}_q} \rho_{q,i}(\boldsymbol{\emptyset}) x_{q,i} \tag{9}$$

*where* $\rho_{q,i}(\boldsymbol{\emptyset})$ *is the marginal gain of adding item* $i$ *to* $q^{th}$ *empty set in* $\boldsymbol{\emptyset}$*, where* $\boldsymbol{\emptyset}$ *is the tuple consisting of* $k$ *empty sets.*

**Theorem 2.** *Given an attacker's solution* $\hat{\mathbf{x}}$ *and associated defender's solution in Formulation (4), denoted by* $\hat{\mathbf{S}} = (\hat{S}_1, \ldots, \hat{S}_k)$ *with an arbitrary permutation of elements in* $\hat{S}_q := \{i_{q,1}, \ldots, i_{q,T_q}\}$ *for* $q \in \{1, \ldots, k\}$*. The following inequality is valid for problem (4) for any* $\mathbf{x} \in \mathcal{X}$ *and it dominates inequality (9), i.e., it provides a tighter lower bound approximation in comparison to (9):*

$$\Phi_D(\mathbf{x}) \geq \Phi_D(\hat{\mathbf{x}}) - \sum_{q=1}^{k} \sum_{t=1}^{T_q} \rho_{q,i_{q,t}}(\hat{\mathbf{S}}_{q,(t)}) x_{q,i_{q,t}}, \tag{10}$$

*where* $\hat{\mathbf{S}}_{q,(t)} = (\hat{S}_1, \ldots, \hat{S}'_{q,(t)}, \emptyset, \ldots, \emptyset)$ *such that* $\hat{S}'_{q,(t)} = \{i_{q,1}, \ldots, i_{q,t-1}\}$ *for* $2 \leq t \leq T_q$ *and* $\hat{\mathbf{S}}_{1,(1)} = \boldsymbol{\emptyset}$*.*

**Theorem 3.** *Given an attacker's solution* $\hat{\mathbf{x}}$ *and associated defender's optimal solutions* $\{\hat{\mathbf{S}}^\omega\}_{\omega \in \Omega}$ *in formulation (6). The following inequality is valid for all* $\mathbf{x} \in \mathcal{X}$*:*

$$\Phi_N(\mathbf{x}) \geq \Phi_N(\hat{\mathbf{x}}) - \sum_{q=1}^{k} \sum_{\omega \in \Omega} \sum_{t=1}^{|\hat{S}_q^\omega|} \overline{p}_\omega \rho^\omega_{q,i^\omega_{q,t}}(\hat{\mathbf{S}}^\omega_{q,(t)}) \xi^\omega_{i^\omega_{q,t}} x_{q,i^\omega_{q,t}}, \tag{11}$$

*where* $\{\overline{p}_\omega\}_\omega \in \mathfrak{P}$ *in the case where* $|\mathfrak{P}| = 1$*, and* $\hat{\mathbf{S}}^\omega_{q,(t)}$ *is defined similar to* $\hat{\mathbf{S}}_{q,(t)}$ *in Theorem 2.*

**Theorem 4.** *Given an attacker's solution* $\hat{\mathbf{x}}$ *and associated defender's optimal solution* $\{\hat{\mathbf{S}}^\omega\}_{\omega \in \Omega}$ *in formulation (7), the following inequality is valid for all* $\mathbf{x} \in \mathcal{X}$*:*

$$\Phi_R(\mathbf{x}) \geq \Phi_R(\hat{\mathbf{x}}) - \sum_{q=1}^{k} \sum_{i \in N} \left( \max_{P \in \mathfrak{P}} \sum_{\omega \in \Omega} p_\omega \rho^\omega_{q,i}(\boldsymbol{\emptyset}) \hat{y}^\omega_{q,i} \right) \xi^\omega_i x_{q,i}, \tag{12}$$

*where for* $i \in N$ *and* $q \in \{1, \ldots, k\}$*,*

$$\hat{y}^\omega_{q,i} = \begin{cases} 1 & \text{if } i \in \hat{S}^\omega_q, \\ 0 & \text{otherwise.} \end{cases}$$

**Theorem 5.** *Given an attacker's solution* $\hat{\mathbf{x}}$ *and associated defender's optimal solutions* $\{\hat{\mathbf{S}}^\omega\}_{\omega \in \Omega}$ *in formulation (8). Then, the following inequality is valid for all*

$\mathbf{x} \in \mathcal{X}$*:*

$$\Phi_A(\mathbf{x}) \geq \Phi_A(\hat{\mathbf{x}}) - \sum_{q=1}^{k} \sum_{\omega\in\Omega} \sum_{t=1}^{|\hat{S}_q^\omega|} \hat{p}_\omega \rho^\omega_{q,i^\omega_{q,t}} (\hat{\mathbf{S}}^\omega_{q,(t)}) \xi^\omega_{i^\omega_{q,t}} x_{q,i^\omega_{q,t}}, \tag{13}$$

*where* $\{\hat{p}_\omega\}_{\omega\in\Omega} \in \arg\max_{P\in\mathfrak{P}} \sum_{\omega\in\Omega} p_\omega f^\omega(\hat{\mathbf{S}}^\omega)$ *and* $\hat{\mathbf{S}}^\omega_{q,(t)}$ *is defined similar to* $\hat{\mathbf{S}}_{q,(t)}$ *in Theorem 2.*

**Observation 6.** Inequalities (10), (11), (12), and (13) are tight (hold at equality) for $\mathbf{x} = \hat{\mathbf{x}}$.

## 4.2 Decomposition Algorithms for DRR and DRO $k$-SIPs

Since formulations (4) and (6) are special cases of DRR $k$-SIP, i.e., (7), we present decomposition algorithm for the latter to avoid repetition. Our decomposition algorithm incorporates inequalities (12) to refine the lower bound for problem (7) during each iteration, as detailed in Algorithm 1. The algorithm begins with an initialization stage where the iteration counter $L = 1$, the upper bound $\theta^{ub}_{drr} \leftarrow \infty$, the lower bound $\theta^{lb}_{drr} \leftarrow -\infty$, and an initial feasible solution $\hat{\mathbf{x}}^1$ is selected from the set $\mathcal{X}$. In each iteration $L \geq 1$, we set $\hat{\mathbf{x}}$ to $\hat{\mathbf{x}}^L$ and solve the defender's problem for every scenario $\omega \in \Omega$, as outlined in Line 4, to obtain $f^\omega(\hat{\mathbf{S}}^\omega)$. Subsequently, in lines 6 and 7, we obtain an extremal optimal probability distribution by solving distribution separation problem:

$$\{\hat{p}_\omega\}_{\omega\in\Omega} \in \arg\min_{P\in\mathfrak{P}} \sum_{\omega\in\Omega} p_\omega f^\omega(\hat{\mathbf{S}}^\omega)$$

and then compute $\Phi_R(\hat{\mathbf{x}}) = \sum_{\omega\in\Omega} \hat{p}_\omega f^\omega(\hat{\mathbf{S}}^\omega)$. Notice that the distribution separation problem is a linear program when the ambiguity set is a polyhedron. In the literature, it has been showcased that many well-known families of ambiguity sets such as moment-based, Wasserstein distance-based, and total variation distance-based can be characterized using linear constraints.

If $\Phi_R(\hat{\mathbf{x}})$ is smaller than the current upper bound $\theta^{ub}_{drr}$, we update $\theta^{ub}_{drr}$ to the value of $\Phi_R(\hat{\mathbf{x}})$ and the best known solution $\hat{\mathbf{x}}^*$ to $\hat{\mathbf{x}}$. Subsequently, we derive and add inequality (12) to a lower bound approximation model from the previous iteration, $\mathcal{M}^{L-1}_{DRR}$, obtaining $\mathcal{M}^L_{DRR}$:

$$\theta^{lb}_{drr} := \min_{\mathbf{x}\in\mathcal{X}} \eta \tag{14}$$

$$\text{s.t. } \eta \geq \Phi_R(\hat{\mathbf{x}}) - \sum_{q=1}^{k} \sum_{i\in N} \left( \max_{P\in\mathfrak{P}} \sum_{\omega\in\Omega} p_\omega \rho^\omega_{q,i}(\emptyset) \hat{y}^\omega_{q,i} \right) \xi^\omega_i x_{q,i} \tag{15}$$

$$\text{for } \hat{\mathbf{x}} \in \{\hat{\mathbf{x}}^1, \dots, \hat{\mathbf{x}}^L\}$$

which is a tighter lower bound approximation. Note that to derive inequality in Line 11, we solve distribution separation problem for $i \in N$ and $q \in \{1, \dots, k\}$. Then, we solve $\mathcal{M}^L_{DRR}$, $L \geq 1$, in Line 12 to get an optimal solution $(\eta^{L+1}, \hat{\mathbf{x}}^{L+1})$ and update the best-known lower bound $\theta^{lb}_{drr}$ to $\eta^{L+1}$. We terminate the algorithm when the optimality gap $(\theta^{ub}_{drr} - \theta^{lb}_{drr})$ lies within a predetermined threshold $\epsilon$.

---

**Algorithm 1** Decomposition Method for the DRR $k$-SIP (7)

---

1: Let $L \leftarrow 1$, $\theta^{lb}_{drr} \leftarrow -\infty$, $\theta^{ub}_{drr} \leftarrow \infty$, $\hat{\mathbf{x}} \leftarrow \hat{\mathbf{x}}^1 \in \mathcal{X}$;
2: **while** $\theta^{ub}_{drr} - \theta^{lb}_{drr} > \epsilon$ **do**
3:     **for** $\omega \in \Omega$ **do**
4:         Solve subproblem to get $f^\omega(\hat{\mathbf{S}}^\omega)$;
5:     **end for**
6:     Compute $\{\hat{p}\}_{\omega \in \Omega} \in \arg\min_{P \in \mathfrak{P}} \sum_{\omega \in \Omega} p_\omega f^\omega(\hat{\mathbf{S}}^\omega)$;
7:     Obtain $\Phi_R(\hat{\mathbf{x}}) = \sum_{\omega \in \Omega} \hat{p}_\omega f^\omega(\hat{\mathbf{S}}^\omega)$ ;
8:     **if** $\theta^{ub}_{drr} > \Phi_R(\hat{\mathbf{x}})$ **then**
9:         $\theta^{ub}_{drr} \leftarrow \Phi_R(\hat{\mathbf{x}})$ and $\hat{\mathbf{x}}^* \leftarrow \hat{\mathbf{x}}$;
10:     **end if**
11:     Add the following inequality in $\mathcal{M}^{L-1}_{DRR}$ to get $\mathcal{M}^L_{DRR}$ :

$$\eta \geq \Phi_R(\hat{\mathbf{x}}) - \sum_{q=1}^{k} \sum_{i \in N} \left( \max_{P \in \mathfrak{P}} \sum_{\omega \in \Omega} p_\omega \rho^\omega_{q,i}(\mathbf{0}) \hat{y}^\omega_{q,i} \right) \xi^\omega_i x_{q,i};$$

12:     Solve $\mathcal{M}^L_{DRR}$ to get optimal solution $(\eta^{L+1}, \hat{\mathbf{x}}^{L+1})$;
13:     Update the lower bound $\theta^{lb}_{drr} \leftarrow \eta^{L+1}$ and $\hat{\mathbf{x}} \leftarrow \hat{\mathbf{x}}^{L+1}$;
14:     $L \leftarrow L + 1$;
15: **end while**
16: **Return**: Optimal Solution Value $\theta^{ub}_{drr}$ and Optimal Solution $\hat{\mathbf{x}}^*$.

---

**Theorem 7.** *Algorithm 1 solves DRR $k$-SIP to global optimality in finite iterations if distribution separation problem, line 6 in Algorithm 1, can be solved in finite iterations.*

*Proof.* Refer to Section 4.3. □

**Remark 2.** The assumption in Theorem 7 holds when the ambiguity set is either polyhedral or defined by mixed-binary constraints, as a finitely converging oracle exists for solving the corresponding distribution separation problem.

**Remark 3. Algorithms for DRO, Deterministic, and Risk-Neutral $k$-SIP.** Algorithm 1 with inequalities (10),(11) and (13) in Line 11 provide an exact algorithm for solving problem (4), (6) and (8) respectively.

## 4.3 Proof of Theorems 1 - 5 and Theorem 7

**Lemma 8** ([29])**.** Let $f$ be a monotone $k$-submodular function. For any $\mathbf{Y}, \mathbf{Z} \in \mathbb{X}(N, k)$, the following inequality holds:

$$f(\mathbf{Y}) \leq f(\mathbf{Z}) + \sum_{q=1}^{k} \sum_{i \in Y_q \setminus \bigcup_{r=1}^{k} Z_r} \rho_{q,i}(\mathbf{Z}) + \sum_{q=1}^{k} \sum_{p \in \{1,\ldots,k\} \setminus \{q\}} \sum_{i \in Y_q \cap Z_p} \rho_{q,i}(\boldsymbol{\emptyset}). \tag{16}$$

**Proof of Theorem 1.** For a feasible attacker's solution $\mathbf{x} = (x_1, \ldots, x_k)$, let $N_{\mathbf{x}} = \cup_{q=1}^{k} N_{x_q}$ be the set of items interdicted by $\mathbf{x}$ where $N_{x_q} = \{i \in N : x_{q,i} = 1\}$ for $q \in \{1, \ldots, k\}$. We define $\mathbf{S} = (S_1, \ldots, S_k)$, where $S_q = \hat{S}_q \setminus N_{x_q}$ for $q \in \{1, \ldots, k\}$. From the definition, it follows that $S_q$ is a subset of $\hat{S}_q$, and $\sum_{i=1}^{n} s_{q,i} \leq \sum_{i=1}^{n} \hat{s}_{q,i} \leq D_q$ holds. Therefore, $\mathbf{S}$ is a feasible defender's solution for the attacker's solution $\mathbf{x}$. Subsequently, referring to (16), we derive the following:

$$\begin{aligned}
f(\hat{\mathbf{S}}) &\leq f(\mathbf{S}) + \sum_{q=1}^{k} \sum_{i \in \hat{S}_q \setminus \bigcup_{r=1}^{k} S_r} \rho_{q,i}(\mathbf{S}) + \sum_{q=1}^{k} \sum_{p \in \{1,\ldots,k\} \setminus \{q\}} \sum_{i \in \hat{S}_q \cap S_p} \rho_{q,i}(\boldsymbol{\emptyset}) && (17)\\
&= f(\mathbf{S}) + \sum_{q=1}^{k} \sum_{i \in \hat{S}_q \setminus \bigcup_{r=1}^{k} S_r} \rho_{q,i}(\mathbf{S}) && (18)
\end{aligned}$$

$$\begin{aligned}
&= f(\mathbf{S}) + \sum_{q=1}^{k} \sum_{i \in \hat{S}_q} \rho_{q,i}(\mathbf{S}) x_{q,i} && (19)\\
&\leq f(\mathbf{S}) + \sum_{q=1}^{k} \sum_{i \in \hat{S}_q} \rho_{q,i}(\boldsymbol{\emptyset}) x_{q,i}. && (20)
\end{aligned}$$

For any $p, q \in \{1, \ldots, k\}$ with $p \neq q$, $\hat{S}_q \cap \hat{S}_p = \emptyset$, and since $S_p \subseteq \hat{S}_p$, equality (18) holds. Similarly, we know that $\hat{S}_q \setminus \bigcup_{r=1}^{k} S_r = \hat{S}_q \setminus S_q$ and for $i \in S_q$, $x_{q,i} = 0$, thereby (19) holds. Finally, inequality (20) is established from the diminishing return property of the monotone $k$-submodular function. Thus, we have

$$\Phi_D(\mathbf{x}) \geq f(\mathbf{S}) \geq f(\hat{\mathbf{S}}) - \sum_{q=1}^{k} \sum_{i \in \hat{S}_q} \rho_{q,i}(\boldsymbol{\emptyset}) x_{q,i} = \Phi(\hat{\mathbf{x}}) - \sum_{q=1}^{k} \sum_{i \in \hat{S}_q} \rho_{q,i}(\boldsymbol{\emptyset}) x_{q,i}. \tag{21}$$

□

**Proof of Theorem 2.** For an attacker's solution $\mathbf{x}$, let $S_q = \hat{S}_q \setminus N_{x_q}$ and $\mathbf{S}_{q,(t)} = (S_1, \ldots, S'_{q,(t)}, \emptyset, \ldots, \emptyset)$ where $S'_{q,(t)} = \hat{S}'_{q,(t)} \setminus N_{x_q}$ for $2 \leq t \leq T_q$ and $q \in \{1, \ldots, k\}$. From the feasibility of $\mathbf{S}$ for $\mathbf{x}$, the following holds:

$$\begin{aligned}
\Phi_D(\mathbf{x}) \geq f(\mathbf{S}) = &\rho_{1,i_{1,1}}(\mathbf{S}_{1,(1)})(1 - x_{1,i_{1,1}}) + \cdots + \rho_{1,i_{1,T_1}}(\mathbf{S}_{1,(T_1)})(1 - x_{1,i_{1,T_1}})\\
&+ \rho_{2,i_{2,1}}(\mathbf{S}_{2,(1)})(1 - x_{2,i_{2,1}}) + \cdots + \rho_{2,i_{1,T_2}}(\mathbf{S}_{2,(T_2)})(1 - x_{2,i_{2,T_2}})
\end{aligned}$$

$$
\begin{aligned}
&\vdots \\
&+\rho_{k,i_{k,1}}(\mathbf{S}_{k,(1)})(1-x_{k,i_{k,1}})+\cdots+\rho_{k,i_{1,T_k}}(\mathbf{S}_{k,(T_k)})(1-x_{k,i_{k,T_k}}) \\
=&\sum_{q=1}^{k}\sum_{t=1}^{T_q}\rho_{q,i_{q,t}}(\mathbf{S}_{q,(t)})(1-x_{q,i_{q,t}}) \\
\geq&\sum_{q=1}^{k}\sum_{t=1}^{T_q}\rho_{q,i_{q,t}}(\hat{\mathbf{S}}_{q,(t)})(1-x_{q,i_{q,t}}) \qquad (22)\\
=&\,\Phi_D(\hat{\mathbf{x}})-\sum_{q=1}^{k}\sum_{t=1}^{T_q}\rho_{q,i_{q,t}}(\hat{\mathbf{S}}_{q,(t)})x_{q,i_{q,t}}. \qquad (23)
\end{aligned}
$$

From the definition of $\mathbf{S}_{q,(t)}$, $S_q \subseteq \hat{S}_q$ and $S^{'}_{q,(t)} \subseteq \hat{S}^{'}_{q,(t)}$ holds for all $q \in \{1,\ldots,k\}$ and $t \in \{1,\ldots,T_q\}$ and thereby, inequality (22) holds from the diminishing return property of $k$-submodular function. It is worth to note that this inequality dominates inequality (9) since $\rho_{q,i_{q,t}}(\boldsymbol{\emptyset}) \geq \rho_{q,i_{q,t}}(\hat{\mathbf{S}}_{q,(t)})$ for all $t \in \{1,\ldots,T_q\}$ for $q \in \{1,\ldots,k\}$. $\square$

**Proof of Theorem 3.** The proof follows the same arguments as those of Theorem 4 and Theorem 5 with a singleton ambiguity set, i.e., $|\mathfrak{P}| = 1$. $\square$

**Proof of Theorem 4.** Given $\tilde{\mathbf{x}} \in \mathcal{X}$, let $\{\tilde{\mathbf{S}}^{\omega}\}_{\omega\in\Omega}$ be an optimal defender's solution, and $\{\tilde{p}_\omega\}_{\omega\in\Omega}$ be an associated extremal probability distribution, i.e.,

$$
\{\tilde{p}_\omega\}_{\omega\in\Omega} \in \arg\min_{P\in\mathfrak{P}} \sum_{\omega\in\Omega} p_\omega f^{\omega}(\tilde{\mathbf{S}}^{\omega}).
$$

Additionally, for $\omega \in \Omega$, define $N^{\omega}_{\tilde{\mathbf{x}}} = \cup_{q=1}^{k} N^{\omega}_{\tilde{x}_q}$ as the set of items interdicted by the attacker's solution $\tilde{\mathbf{x}}$, where $N^{\omega}_{\tilde{x}_q} = \{i \in N : \tilde{x}_{q,i} \cdot \xi^{\omega}_i = 1\}$ for $q \in \{1,\ldots,k\}$. For $\omega \in \Omega$, we define $\{\overline{\mathbf{S}}^{\omega}\}_{\omega\in\Omega} = (\overline{S}^{\omega}_1,\ldots,\overline{S}^{\omega}_q)$ where $\overline{S}^{\omega}_q = \tilde{S}^{\omega}_q \setminus N^{\omega}_{\tilde{x}_q}$ for all $q \in \{1,\ldots,k\}$, $\omega \in \Omega$. This ensures that $\{\overline{\mathbf{S}}^{\omega}\}_{\omega\in\Omega}$ is a feasible defender's solution for the given $\tilde{\mathbf{x}}$. Then,

$$
\begin{aligned}
\Phi_R(\tilde{\mathbf{x}}) =& \sum_{\omega\in\Omega}\tilde{p}_\omega f^{\omega}(\tilde{\mathbf{S}}^{\omega}) \\
\geq& \sum_{\omega\in\Omega}\tilde{p}_\omega f^{\omega}(\overline{\mathbf{S}}^{\omega}) \qquad (24a)\\
\geq& \sum_{\omega\in\Omega}\tilde{p}_\omega \left( f^{\omega}(\hat{\mathbf{S}}^{\omega}) - \sum_{q=1}^{k}\sum_{i\in\hat{S}^{\omega}_q}\rho^{\omega}_{q,i}(\boldsymbol{\emptyset})\xi^{\omega}_i\tilde{x}_{q,i}\right) \qquad (24b)\\
\geq& \,\Phi_R(\hat{\mathbf{x}}) - \sum_{q=1}^{k}\sum_{\omega\in\Omega}\sum_{i\in\hat{S}^{\omega}_q}\tilde{p}_\omega\rho^{\omega}_{q,i}(\boldsymbol{\emptyset})\xi^{\omega}_i\tilde{x}_{q,i} \qquad (24c)\\
=& \,\Phi_R(\hat{\mathbf{x}}) - \sum_{q=1}^{k}\sum_{i\in N}\sum_{\omega\in\Omega}\tilde{p}_\omega\rho^{\omega}_{q,i}(\boldsymbol{\emptyset})\hat{y}^{\omega}_{q,i}\xi^{\omega}_i\tilde{x}_{q,i} \qquad (24d)
\end{aligned}
$$

$$\geq \Phi_R(\hat{\mathbf{x}}) - \sum_{q=1}^{k} \sum_{i \in N} \left( \max_{P \in \mathfrak{P}} \sum_{\omega \in \Omega} p_\omega \rho^\omega_{q,i}(\emptyset) \hat{y}^\omega_{q,i} \right) \xi^\omega_i \tilde{x}_{q,i}. \tag{24e}$$

Inequality (24a) holds from the feasibility of set $\{\overline{\mathbf{S}}^\omega\}_{\omega\in\Omega}$ for the attacker's solution $\tilde{\mathbf{x}}$ and inequality (24b) follows from (25). Since $\Phi_R(\hat{\mathbf{x}}) = \min_{P\in\mathfrak{P}} \mathbb{E}_P[Q_\omega(\hat{\mathbf{x}}, \xi^\omega)]$, inequality (24c) holds. Also, from the definition of $\hat{y}^\omega_{q,i}$, inequality (21) holds. Finally, inequality (24e) holds because $\rho^\omega_{q,i}(\emptyset), \hat{y}^\omega_{q,i}, \xi^\omega_i$ and $\tilde{x}_{q,i}$ are non-negative for all $q \in \{1,\ldots,k\}$ and $i \in N$ and $\{\tilde{p}_\omega\}_{\omega\in\Omega}$ is a feasible solution for the maximization problem in inequality (24e). $\square$

**Proof of Theorem 5.** For any feasible attacker's solution $\mathbf{x}$, define $\overline{\mathbf{S}}^\omega = (\overline{S}^\omega_1, \ldots, \overline{S}^\omega_q)$ where $\overline{S}^\omega_q = \hat{S}^\omega_q \setminus N^\omega_{x_q}$ for all $q \in \{1,\ldots,k\}$ and $\omega \in \Omega$. From the definition, $\{\overline{\mathbf{S}}^\omega\}_{\omega\in\Omega}$ is a feasible defender's solution for the given $\mathbf{x}$. Referring to (21), the following holds for all $\omega \in \Omega$:

$$f^\omega(\overline{\mathbf{S}}^\omega) \geq f^\omega(\hat{\mathbf{S}}^\omega) - \sum_{q=1}^{k} \sum_{i \in \hat{S}^\omega_q} \rho^\omega_{q,i}(\emptyset) \xi^\omega_i x_{q,i}. \tag{25}$$

Moreover, from (22), the following inequality holds for all $\omega \in \Omega$:

$$f^\omega(\overline{\mathbf{S}}^\omega) \geq f^\omega(\hat{\mathbf{S}}^\omega) - \sum_{q=1}^{k} \sum_{t=1}^{T^\omega_q} \rho^\omega_{q,i^\omega_{q,t}}(\hat{\mathbf{S}}^\omega_{q,(t)}) \xi^\omega_{i^\omega_{q,t}} x_{q,i^\omega_{q,t}}. \tag{26}$$

Let $\{\bar{p}_\omega\}_{\omega\in\Omega} \in \arg\max_{P\in\mathfrak{P}} \sum_{\omega\in\Omega} p_\omega f(\overline{\mathbf{S}}^\omega)$. Then,

$$\Phi_A(\mathbf{x}) \geq \sum_{\omega\in\Omega} \bar{p}_\omega f^\omega(\overline{\mathbf{S}}^\omega) \tag{27a}$$

$$\geq \sum_{\omega\in\Omega} \hat{p}_\omega f^\omega(\overline{\mathbf{S}}^\omega) \tag{27b}$$

$$\geq \sum_{\omega\in\Omega} \hat{p}_\omega \left( f^\omega(\hat{\mathbf{S}}^\omega) - \sum_{q=1}^{k} \sum_{t=1}^{|\hat{S}^\omega_q|} \rho^\omega_{q,i^\omega_{q,t}}(\hat{\mathbf{S}}^\omega_{q,(t)}) \xi^\omega_{i^\omega_{q,t}} x_{q,i^\omega_{q,t}} \right) \tag{27c}$$

$$= \Phi_A(\hat{\mathbf{x}}) - \sum_{q=1}^{k} \sum_{\omega\in\Omega} \sum_{t=1}^{|\hat{S}^\omega_q|} \hat{p}_\omega \rho^\omega_{q,i^\omega_{q,t}}(\hat{\mathbf{S}}^\omega_{q,(t)}) \xi^\omega_{i^\omega_{q,t}} x_{q,i^\omega_{q,t}}. \tag{27d}$$

From the feasibility of $\overline{\mathbf{S}}^\omega$ for given $\mathbf{x}$, (27a) holds. Inequality (27b) is based on the definition of $\{\bar{p}_\omega\}_{\omega\in\Omega}$. From (26), (27c) and (27d) follow. $\square$

**Proof of Theorem 7.** Given that $\Omega$ is a set of finite possible realizations, we solve finite $|\Omega|$ number of subproblems in each iteration of the while loop (steps from line 2 to line 15). For each subproblem (a $k$-submodular maximization problem), an optimal solution can be found in a finite iterations, given that the ground set $N = \{1,\ldots,n\}$ has a finite number of elements [29]. Additionally, we solve the master problem (a mixed-binary problem) once and solve distribution separation problem in line 6 and 7) at most $nk$ number of times. Assuming the existence of an algorithm that converges in finite iterations for distribution separation problem (line 6 and 7), the steps from line 2 to line 15 will also take finite number of iterations for any counter $L$. Now, we only need to prove that steps from line 2 to line 15 iterate finite number of times, i.e., $L$ is finite, until the algorithm converges to an optimal solution. Consider an iteration $L = t < \infty$. Suppose we obtain optimal solution $(\mathbf{x}^t, \eta^t)$ by solving the master problem $\mathcal{M}^t_{DRR}$, where $\eta^t = \theta^{lb}_{drr}$. If there is no optimality cut separating $(\mathbf{x}^t, \eta^t)$, then this solution is feasible for the DRO $k$-SIP, implying $\eta^t = \Phi_R(\mathbf{x}^t) = \theta^{ub}_{drr}$. Hence, $\eta^t = \theta^{lb}_{drr} = \theta^{ub}_{drr}$, which satisfies termination condition, and we have optimal solution in the finite iterations (at iteration $t$).

If there exists an optimality cut separating $(\mathbf{x}^t, \eta^t)$, it will be added as in line 11 and the process (while loop) continues. Assume we obtain an optimal master solution $(\mathbf{x}^u, \eta^u)$ in iteration $u$ where $u > t$. If we already visit $\mathbf{x}^u$ in previous iterations, this implies $(\mathbf{x}^u, \eta^u)$ is a feasible solution of the DRR $k$-SIP and no new optimality cut will be added. We already show that in this case, algorithm will be terminated with the optimal solution $(\mathbf{x}^u, \eta^u)$. Finally, if $(\mathbf{x}^u, \eta^u)$ is the new solution we never visited before, the new optimality cut will be added and algorithm will be continued. Since the attacker's solution $\mathbf{x}$ is a vector of binary variables, there exists finite number of feasible attacker's solution that can be visited. Consequently, the algorithm will be terminated in finite iterations with an optimal solution because of Observation 6. $\square$

## 5 Results of Computational Experiments

In this section, we present the computational results for solving DRO $k$-SIP, DRR $k$-SIP, risk-neutral (stochastic) $k$-SIP, and deterministic $k$-SIP instances. The evaluation is divided into two categories based on the objective function of the defender: $k = 1$ and $k = 2$ submodular function in feature selection interdiction problem (Section 5.1) and weighted coverage interdiction problem (Section 5.2), respectively. For both problems, we set $B_i = 1$ for all $i$ in the constraint 4b, thereby limiting the attacker to interdicting at most one feature or sensor type per element of the ground set. The algorithms were implemented in Python 3.8.5 with the Gurobi 9.5.2 optimization solver and the experiments were conducted using Intel Xeon(R) W-2255 processor (3.7 GHz) with 32GB RAM. We consider two types of ambiguity set that are widely used in the literature: Moment matching set and Wasserstein ambiguity set.

The moment matching set, denoted by $\mathfrak{P}_M$, is defined by constraining the first moments of a family of $m$ predefined functions, $f_\lambda : \Omega \to \mathbb{R}$, for $\lambda = 1, \ldots, m$. These functions are considered real-valued and measurable with respect to the sample space $\Omega$ and its sigma-algebra $\mathcal{F}$ (on which the expectation operator $\mathbb{E}[\cdot]$ is defined). For each scenario $\omega \in \Omega$, $f_\lambda(\omega)$ represents the realized value of the $\lambda$-th random quantity of interest (e.g., $f_\lambda(\omega)$ could represent elements such as $\xi^\omega_i$ from the success of attack vector or specific feature values derived from the realized dataset $\mathcal{D}^\omega$.

For each such function $f_\lambda$ (which represents a random variable), its empirical mean (i.e., sample average) can be calculated as $\overline{f}_\lambda = \sum_{\omega \in \Omega} f_\lambda(\omega)/|\Omega|$. Given a tolerance parameter $\epsilon_M \geq 0$ for each $f_\lambda$, corresponding lower and upper bounds for its true first moment, denoted $\underline{u}_\lambda$ and $\bar{u}_\lambda$ respectively, can, for example, be obtained as $\underline{u}_\lambda = (1 - \epsilon_M)\overline{f}_\lambda$ and $\bar{u}_\lambda = (1 + \epsilon_M)\overline{f}_\lambda$. The moment matching ambiguity set $\mathfrak{P}_M$ is then defined as:

$$\mathfrak{P}_M := \left\{ \{p_\omega\}_{\omega\in\Omega} : \underline{u}_\lambda \leq \mathbb{E}_P[f_\lambda] \leq \bar{u}_\lambda \quad \forall \lambda = 1, \ldots, m, \sum_{\omega\in\Omega} p_\omega = 1; p_\omega \geq 0, \ \forall \omega \in \Omega \right\}, \tag{28}$$

where $\mathbb{E}_P[f_\lambda] = \sum_{\omega\in\Omega} p_\omega f_\lambda(\omega)$ represents the expected value of the random variable $f_\lambda$ under a specific probability distribution $P = \{p_\omega\}_{\omega\in\Omega}$.

In contrast, the Wasserstein ambiguity set is a set of probability distributions within a given $\epsilon_W (> 0)$ distance from a reference distribution $P^* = \{p^*_\omega\}_{\omega\in\Omega}$ and is defined as:

$$\begin{aligned} \mathfrak{P}_W := \Bigg\{ P = \{p_\omega\}_{\omega\in\Omega} : \sum_{\omega\in\Omega} p_\omega = 1; \sum_{\omega_i \neq \omega_j \in\Omega} \|\omega_i - \omega_j\|_1 v_{\omega_i,\omega_j} \leq \epsilon_W; \\ \sum_{\omega_j\in\Omega} v_{\omega_i,\omega_j} = p_{\omega_i}, \ \text{ for all } \omega_i \in \Omega; \sum_{\omega_i\in\Omega} v_{\omega_i,\omega_j} = p^*_{\omega_j}, \forall \omega_j \in \Omega; \\ p_\omega \geq 0, \ \text{ for all } \omega \in \Omega; \ v_{\omega_i,\omega_j} \geq 0, \forall \omega_i, \omega_j \in \Omega \Bigg\}. \end{aligned} \tag{29}$$

## 5.1 Feature Selection Interdiction Problem

### 5.1.1 Instance Generation and setting

For FSIP instances, we use Wisconsin Breast Cancer Data set ($\mathcal{D}$), which is widely recognized in the field of medical informatics for research purposes [21]. This data set is available from the UCI Machine Learning Repository [48, 49] and consists of 569 data points each representing a patient. Each patient is described by $n = 30$ features derived from digitized images of breast mass biopsies. The data set is primarily used for binary classification tasks within machine learning frameworks. We assume that each data point inherently possesses uncertainty, indicating potential noise within each data point. Specifically, the “true” data point exists within a ball of radius $\delta (\geq 0)$ that is centered around each data point in $\mathcal{D}$.

Utilizing the data set $\mathcal{D}$, we generated multiple data sets, $\{\mathcal{D}^\omega\}_{\omega\in\Omega}$, to represent a set of distinct possible realizations of the true data. For each data set $\mathcal{D}^\omega$, $\omega \in \Omega$, every feature value of the each observed data in $\mathcal{D}$ was perturbed by adding uniformly distributed random noise within the range defined by $\delta = 0.1$. This means that in $\mathcal{D}^\omega$, $\omega \in \Omega$, all data points are perturbed independently but concurrently based on the original data set $\mathcal{D}$, ensuring that we generated distinct data sets without cross-combining perturbed data points from different realizations. Furthermore, corresponding to each variant of the data set, we generate sets of binary vectors, $\xi^\omega \in \{0,1\}^n$, $\omega \in \Omega$, indicating the success of the attack on feature $i$ if $\xi^\omega_i = 1$. These vectors are generated from the Bernoulli distribution with a probability 0.75, meaning $\xi_i^\omega = 1$ with probability 0.75 for $i \in \{1, \ldots, n\}$. Finally, we fix the number of scenarios $|\Omega|$ to 100.

### 5.1.2 In sample test: Impact of Solutions on Accuracy of SVC Models

Using $\{\mathcal{D}^\omega\}_{\omega\in\Omega}$, we solve DRO, Risk-Neutral and DRR 1-SIP using the proposed decomposition algorithms and obtain optimal solutions for attackers ($\mathbf{x}_{RA}$, $\mathbf{x}_{RN}, \mathbf{x}_{RR}$) and for defenders across different data realizations ($\{\mathbf{S}^\omega_{RA}\}_{\omega\in\Omega}, \{\mathbf{S}^\omega_{RN}\}_{\omega\in\Omega}, \{\mathbf{S}^\omega_{RR}\}_{\omega\in\Omega}$). We also solve the deterministic problem (4) using original data set $\mathcal{D}$ to find an optimal attacker's solution $\mathbf{x}_{DT}$. Subsequently, for $\{\mathcal{D}^\omega\}_{\omega\in\Omega}$, we obtained a set of defender's solution $\{\mathbf{S}^\omega_{DT}\}_{\omega\in\Omega}$ when attacker's solution is fixed to $\mathbf{x}_{DT}$. Then, for each scenario $\omega \in \Omega$, we partitioned the data set $\mathcal{D}^\omega$ into training and test sets using a 70:30 split. We trained four SVCs using the train data from $\mathcal{D}^\omega$, however, each with defender selected features $\mathbf{S}^\omega_{DT}, \mathbf{S}^\omega_{RA}, \mathbf{S}^\omega_{RN}$ and $\mathbf{S}^\omega_{RR}$ respectively. The impact of attacking strategies was assessed by comparing the test accuracy of the SVC models for test data from $\mathcal{D}^\omega$. For the moment matching set with $\epsilon_M = 0.05$, the results are visualized in Figure 2 where each dot represents test accuracy of SVC that is trained with defender selected features under each model setting across all $\omega \in \Omega$. Each violin plot's outline encodes the full distribution of accuracy values, with colors distinguishing the four models (Deterministic, DRO, Risk-neutral, DRR). Inside each plot, the horizontal bar at the center marks the mean accuracy, while the top and bottom ticks indicate the maximum and minimum accuracies, which are also annotated numerically.

In Figure 2(a), observe that the attacker's solution $\mathbf{x}_{RA}$ achieves the lowest maximum test accuracy, 95%, among the four evaluated problems, indicating its robustness from the attacker's perspective, as higher accuracy is undesirable for the attacker. However, this solution could not degrade model performance as significantly as the attacker's solution $\mathbf{x}_{RR}$, which can reduce performance to 82% but also paradoxically permits accuracy up to 96% for some realizations. This result aligns with our expectation that attacker would opt for $\mathbf{x}_{RR}$ (distributionally risk-receptive solution) when they aim to minimize the defender's objective and to reduce test accuracy as much as possible which comes with the risk of allowing the defender to achieve their objective that is not attainable if the attacker opts for any other attacking strategies. In Figure 2(b), we observe outcomes for different budget and notice that by not leveraging the inherent uncertainties present in real-world data sets, $\mathbf{x}_{DT}$ strategy fails to effectively degrade the model's performance in comparison to the other three attacking strategies. This (deterministic) attacking strategy leads to a relatively stable but less impactful outcome from an attacker's perspective while the other three (stochastic) attacking solutions force a wider range of model performance, showing their ability to substantially influence model accuracy or compromise data integrity. Overall, from the defender's perspective, solving the DRR $k$-SIP offers insights by identifying the most vulnerable features, whose removal could significantly compromise data quality.

In Tables 1 and 2, we report the objective function values, $\Phi_A(\mathbf{x}_{RA})$, $\Phi_N(\mathbf{x}_{RN})$, and $\Phi_R(\mathbf{x}_{RR})$, and the total solution time for instances with equal attacker and defender budgets from 1 to 14. Specifically, the solution time is the CPU time required

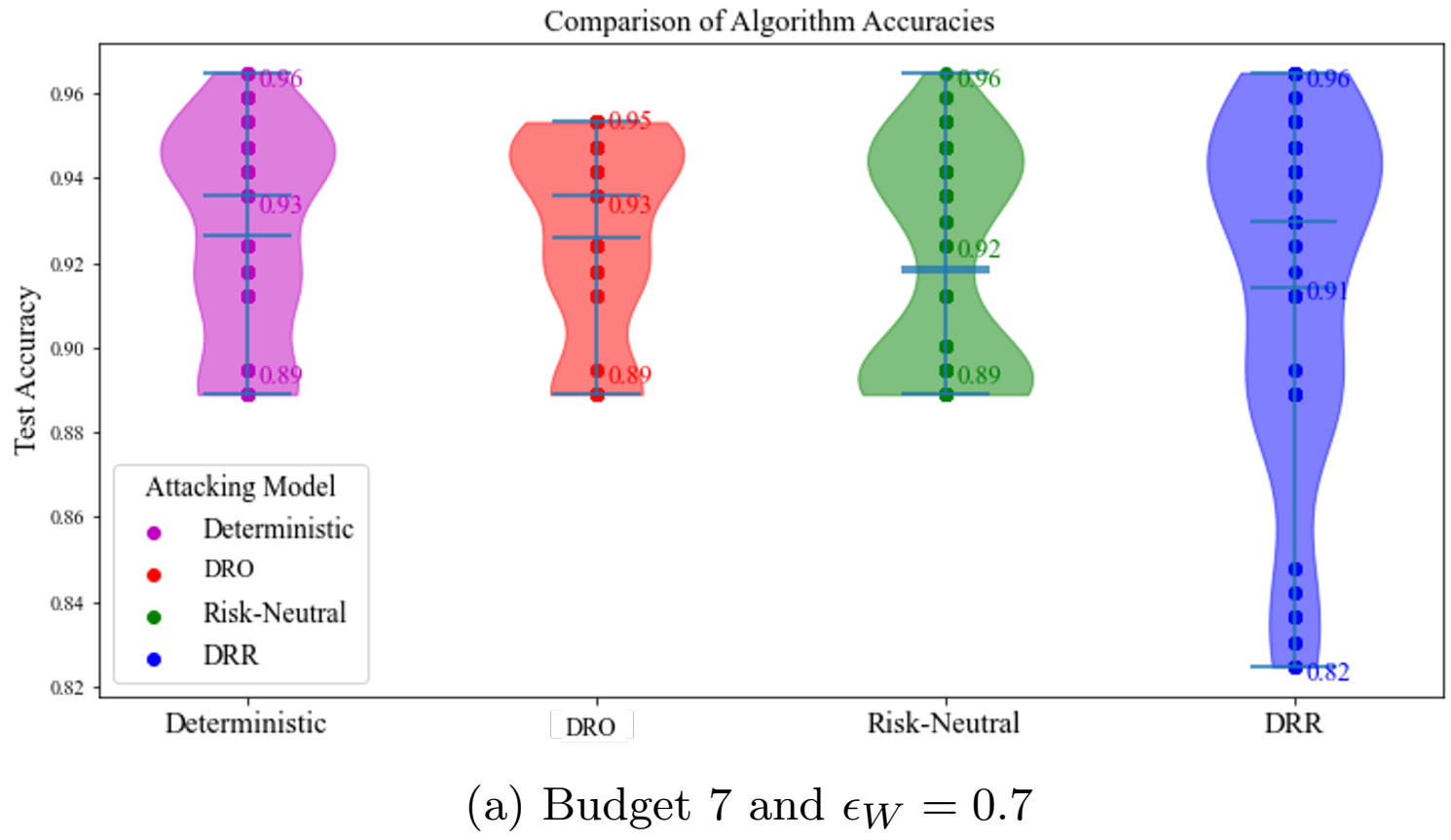

(a) Budget 7 and $\epsilon_W = 0.7$

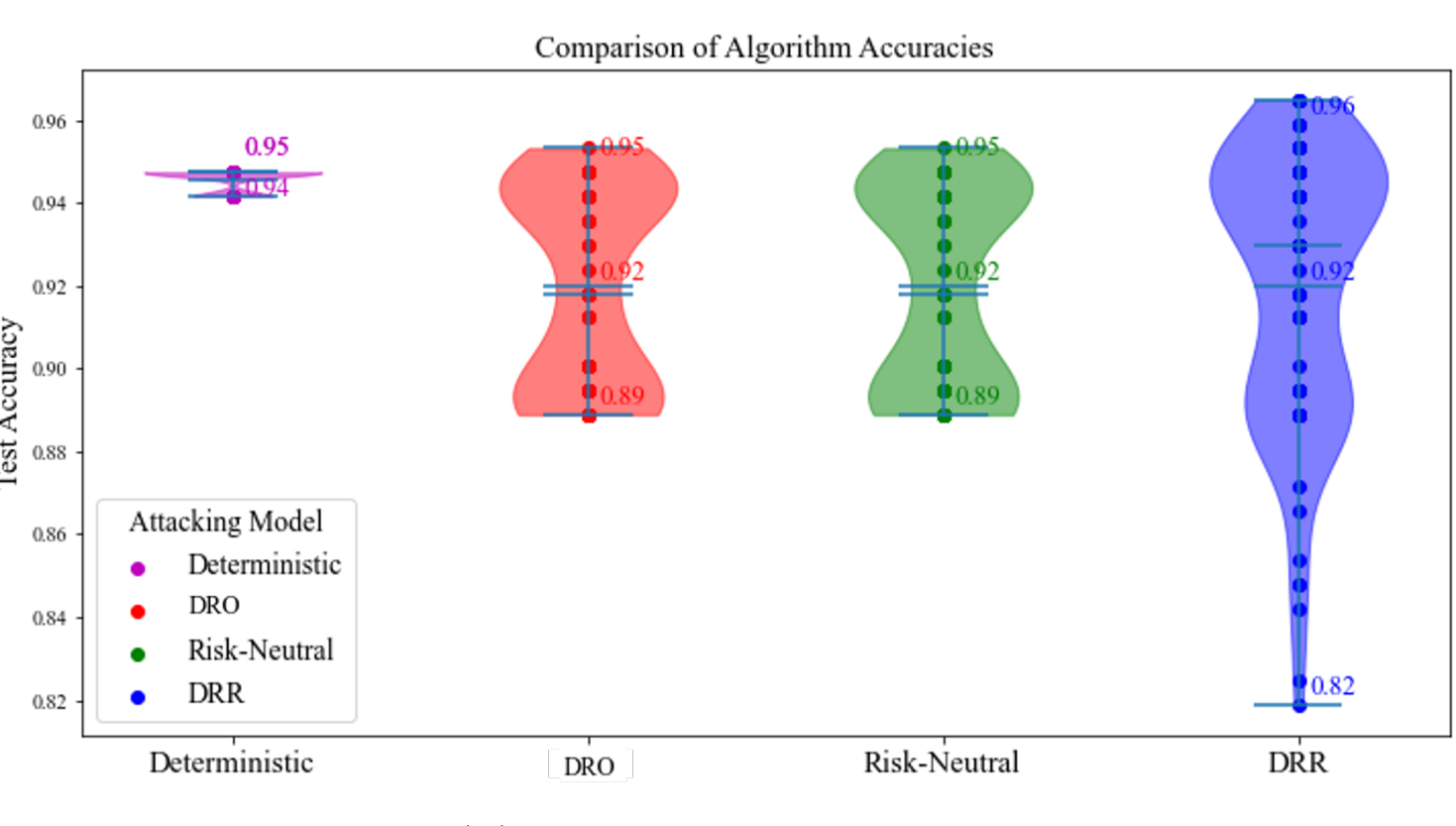

(b) Budget 8 and $\epsilon_W = 0.7$

**Fig. 2**: Performance of Support Vector Classifiers (SVCs) Trained on Defender-Selected Features Across Different Problem Scenarios and Budget Settings

for the algorithm to terminate at an optimal solution (i.e., when the optimality gap $\theta^{ub} - \theta^{lb} \leq \epsilon$) including every step of the algorithm, such as solving the defender's subproblems, solving the distribution separation problem, and deriving the valid inequality. In Figures 3(a) and 3(b), we present the objective function values as we increase $\epsilon_W$ and $\epsilon_M$. Observe that the gap between risk-averse solution values and risk-receptive solution values increase as we increase the $\epsilon_W$ and $\epsilon_M$ to increase the size of the ambiguity set. Also, objective values from DRO $k$-SIP and DRR $k$-SIP establish an interval for the objective value of the Risk-Neutral $k$-SIP, this trend implies that when the attacker has less information associate with the probability distribution, i.e., when they assume the larger ambiguity set, they construct large confidence interval for the expectation of defender's objective value. For the attacker, the interval constructed for the expected values of the defender's objective in a distributional ambiguity, as illustrated in Figure 3(a), serves as a valuable reference which enables them to decide whether to pursue aggressive actions or conservative approach.

| Instance (Budgets) | Wasserstein ambiguity set ($\epsilon_W$) | | | | | | | | | | | |
|---|---|---|---|---|---|---|---|---|---|---|---|---|
| | $\epsilon_W = 0.3$ | | | | | | $\epsilon_W = 0.5$ | | | | | |
| | DRO $k$-SIP | | Risk Neutral $k$-SIP | | DRR $k$-SIP | | DRO $k$-SIP | | Risk Neutral $k$ | | DRR $k$-SIP | |
| | Reward $\Phi_A(\mathbf{x}_{RA})$ | Time (s) | Reward $\Phi_N(\mathbf{x}_{RN})$ | Time (s) | Reward $\Phi_R(\mathbf{x}_{RR})$ | Time (s) | Reward $\Phi_A(\mathbf{x}_{RA})$ | Time (s) | Reward $\Phi_N(\mathbf{x}_{RN})$ | Time (s) | Reward $\Phi_R(\mathbf{x}_{RR})$ | Time (s) |
| 1 | 13.8 | 1.4 | 13.6 | 1.3 | 13.5 | 1.6 | 13.8 | 0.8 | 13.6 | 1.0 | 13.5 | 1.5 |
| 2 | 18.1 | 39.5 | 17.6 | 39.6 | 17.2 | 49.7 | 18.4 | 46.0 | 17.6 | 39.5 | 17.0 | 42.8 |
| 3 | 19.7 | 129.4 | 19.1 | 179.6 | 18.5 | 172.2 | 20.0 | 135.6 | 19.1 | 180.1 | 18.3 | 179.5 |
| 4 | 20.7 | 316.8 | 20.0 | 423.5 | 19.3 | 413.9 | 21.1 | 317.6 | 20.0 | 422.3 | 18.8 | 497.3 |
| 5 | 21.7 | 912.9 | 20.1 | 589.9 | 18.5 | 507.2 | 22.3 | 983.2 | 20.1 | 590.3 | 17.8 | 660.3 |
| 6 | 21.5 | 848.1 | 18.5 | 399.0 | 14.7 | 845.0 | 22.5 | 1799.9 | 18.5 | 400.2 | 12.7 | 599.8 |
| 7 | 22.0 | 2016.7 | 18.6 | 553.9 | 14.3 | 693.7 | 23.0 | 2656.4 | 18.6 | 552.9 | 11.8 | 910.7 |
| 8 | 22.0 | 2074.5 | 18.7 | 1068.6 | 14.4 | 1610.4 | 23.2 | 3451.3 | 18.7 | 1068.2 | 11.9 | 1203.5 |
| 9 | 21.5 | 3376.7 | 17.6 | 1609.1 | 12.3 | 1497.9 | 22.7 | 2639.0 | 17.6 | 1618.6 | 9.2 | 1338.8 |
| 10 | 21.1 | 2701.5 | 17.7 | 1301.2 | 12.4 | 1954.6 | 22.2 | 1788.8 | 17.7 | 1300.7 | 9.4 | 1390.8 |
| 11 | 21.3 | 2782.0 | 16.8 | 1722.2 | 11.4 | 2290.5 | 22.9 | 3779.9 | 16.8 | 1721.8 | 8.3 | 1870.5 |
| 12 | 20.8 | 2684.7 | 17.6 | 2103.6 | 11.9 | 2321.2 | 22.0 | 2696.0 | 17.6 | 2093.0 | 8.6 | 1927.1 |
| 13 | 18.7 | 2246.2 | 16.9 | 1629.1 | 14.9 | 1903.1 | 22.0 | 2077.3 | 16.9 | 1635.8 | 7.9 | 2504.8 |
| 14 | 18.3 | 4220.5 | 16.4 | 2331.7 | 14.5 | 2458.4 | 22.1 | 3937.3 | 16.4 | 2331.8 | 7.7 | 2863.2 |

**Table 1**: Computational results for solving DRO $k$-SIP, Risk-neutral $k$-SIP, and DRR $k$-SIP with Wisconsin Breast Cancer Data with Wasserstein ambiguity set.

| Instance (Budgets ) | Moment matching set ($\epsilon_M$) | | | | | | | | | | | |
|---|---|---|---|---|---|---|---|---|---|---|---|---|
| | $\epsilon_M = 0.03$ | | | | | | $\epsilon_M = 0.05$ | | | | | |
| | DRO $k$-SIP | | Risk Neutral $k$-SIP | | DRR $k$-SIP | | DRO $k$-SIP | | Risk Neutral $k$ | | DRR $k$-SIP | |
| | Reward $\Phi_A(\mathbf{x}_{RA})$ | Time (s) | Reward $\Phi_N(\mathbf{x}_{RN})$ | Time (s) | Reward $\Phi_R(\mathbf{x}_{RR})$ | Time (s) | Reward $\Phi_A(\mathbf{x}_{RA})$ | Time (s) | Reward $\Phi_N(\mathbf{x}_{RN})$ | Time (s) | Reward $\Phi_R(\mathbf{x}_{RR})$ | Time (s) |
| 1 | 13.8 | 3.6 | 13.6 | 1.0 | 13.5 | 54.2 | 13.8 | 3.4 | 13.6 | 1.1 | 13.5 | 51.1 |
| 2 | 17.9 | 36.0 | 17.6 | 39.3 | 17.4 | 115.4 | 18.0 | 35.7 | 17.6 | 38.8 | 17.3 | 122.6 |
| 3 | 19.3 | 158.1 | 19.1 | 179.3 | 18.7 | 304.3 | 19.4 | 156.7 | 19.1 | 178.2 | 18.7 | 298.0 |
| 4 | 20.3 | 364.1 | 20.0 | 425.2 | 19.7 | 775.1 | 20.4 | 386.4 | 20.0 | 419.8 | 19.6 | 750.7 |
| 5 | 20.9 | 795.4 | 20.1 | 593.0 | 19.3 | 802.8 | 21.0 | 928.7 | 20.1 | 587.9 | 19.2 | 812.3 |
| 6 | 20.4 | 927.8 | 18.5 | 407.7 | 15.7 | 624.6 | 20.7 | 1193.4 | 18.5 | 397.3 | 15.1 | 613.2 |
| 7 | 20.6 | 1567.9 | 18.6 | 552.4 | 15.7 | 747.2 | 20.9 | 1528.9 | 18.6 | 552.0 | 15.1 | 657.5 |
| 8 | 20.5 | 1216.4 | 18.7 | 1067.8 | 15.9 | 1264.5 | 20.8 | 2555.8 | 18.7 | 1065.9 | 15.4 | 1286.3 |
| 9 | 19.8 | 2048.0 | 17.6 | 1605.8 | 14.0 | 2031.9 | 20.2 | 2078.5 | 17.6 | 1605.7 | 13.3 | 1656.0 |
| 10 | 19.6 | 1744.1 | 17.7 | 1302.2 | 14.4 | 1953.2 | 19.9 | 1827.8 | 17.7 | 1303.6 | 13.7 | 2029.1 |
| 11 | 19.1 | 2832.8 | 16.8 | 1724.6 | 12.5 | 2273.5 | 19.5 | 3051.6 | 16.8 | 1710.0 | 11.9 | 2586.1 |
| 12 | 19.5 | 2723.9 | 17.6 | 2077.2 | 13.9 | 2353.0 | 19.9 | 3765.9 | 17.6 | 2067.2 | 13.3 | 2373.3 |
| 13 | 18.8 | 2195.8 | 16.9 | 1641.7 | 13.1 | 2736.8 | 19.2 | 2203.9 | 16.9 | 1674.9 | 12.5 | 4221.1 |
| 14 | 18.7 | 2235.5 | 16.4 | 2335.0 | 12.0 | 3864.4 | 19.1 | 1586.1 | 16.4 | 2413.6 | 11.4 | 4650.6 |

**Table 2**: Computational results for solving DRO $k$-SIP, Risk-neutral $k$-SIP and DRR $k$-SIP with Wisconsin Breast Cancer Data with Moment matching set

### 5.1.3 Out of Sample Testing of Deterministic, Risk-Neutral, DRO, and DRR

For out-of-sample testing, we generate another set of possible realizations $\hat{\Omega} := \{\hat{\omega}_1, \ldots, \hat{\omega}_{|\hat{\Omega}|}\}$ of $(\xi, \mathcal{D})$ using Bernoulli and uniform distributions for $\xi$ and $\mathcal{D}$, respectively, as those employed for generating in-sample tests. Then, we mirrored the same procedure used in section 5.1.2 to evaluate the impact of attacking solutions through the test accuracy of SVCs which are trained solely with defender selected features. For these experiments, we increased the number of scenarios to 300, i.e, $|\hat{\Omega}| = 300$. We

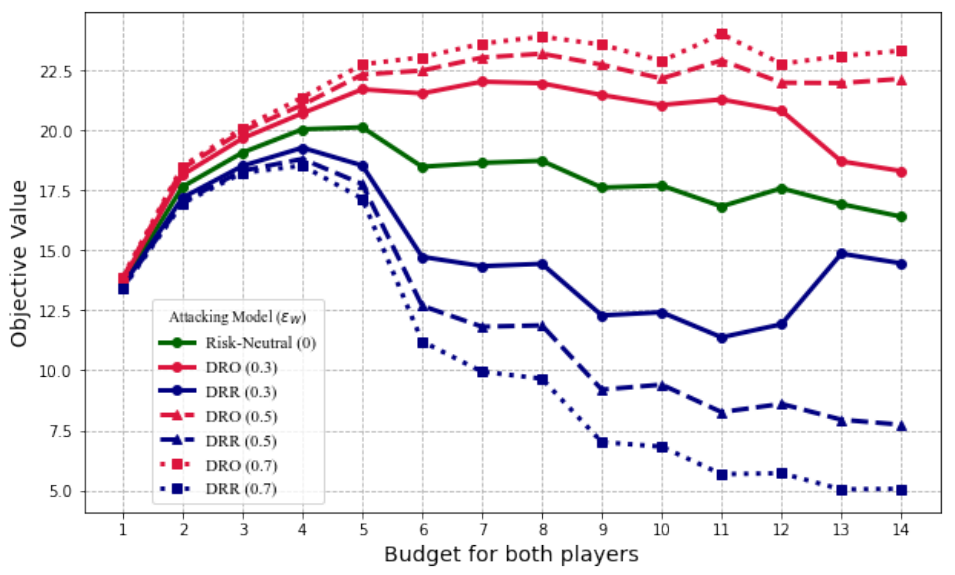

(a) Comparisons of optimal objective values with Wasserstein ambiguity set defined by $\epsilon_W \in \{0.3, 0.5, 0.7\}$

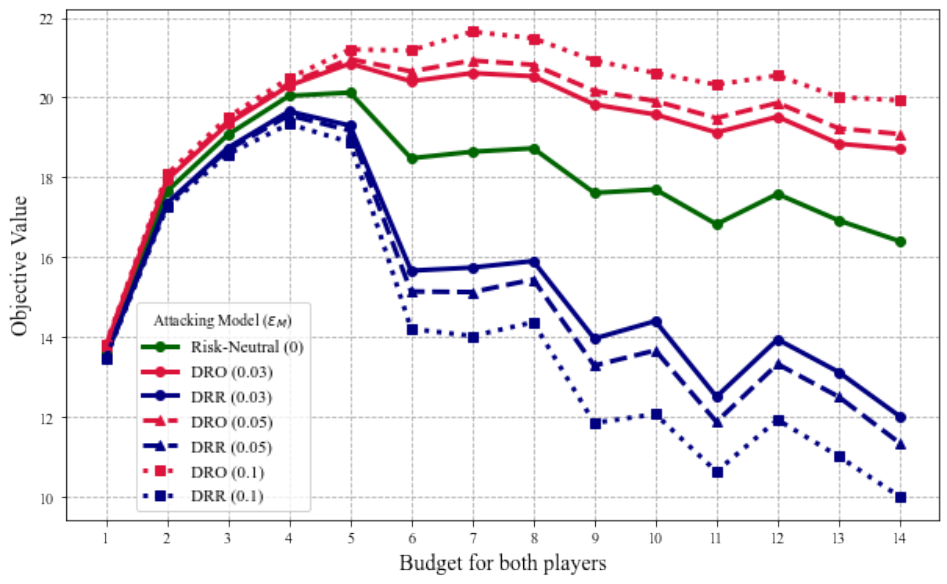

(b) Comparisons of optimal objective values with Moment matching set defined by $\epsilon_M \in \{0.03, 0.05, 0.1\}$

**Fig. 3**: Comparisons of optimal objective values with different ambiguity sets

varied the size of the Wasserstein ambiguity set and for each defined ambiguity set, we followed the aforementioned procedure across different budget settings for both the attacker and the defender. The test results are presented in Figures 4, 5 and 6.

Overall, we notice that the solution obtained using the stochastic models, $\mathbf{x}_{RA}$ and $\mathbf{x}_{RR}$ provide significant advantages over $\mathbf{x}_{DT}$ obtained by solving deterministic model. The $\mathbf{x}_{RA}$ reduces the variability of the outcome, offering a robust attacking strategy that minimizes the risk. Also, Figures 5 and 6 show that $\mathbf{x}_{RR}$ consistently yields lower mean test accuracy. For example, Figure 6(b), the mean accuracy drops from 0.963 under $\mathbf{x}_{DT}$ to 0.953 under $\mathbf{x}_{RR}$. Likewise, at Figure 6(c), the mean accuracy falls from 0.961 to 0.952, confirming that $\mathbf{x}_{RR}$ degrades performance more aggressively than $\mathbf{x}_{DT}$. Theses two distinct strategies offer decision makers the flexibility to choose between robustness and aggressiveness.

***Calibration of the Wasserstein Radius.***

To identify the radius $\epsilon_W$ that maximizes each attacker's out-of-sample impact, we consider the candidate set $\epsilon_W \in \{0.1, 0.3, 0.5, 0.7, 1.0, 1.5\}$. For each $\epsilon_W$ and each attacker model $m \in \{\text{RA}, \text{RN}, \text{RR}\}$ (DRO, risk-neutral, DRR), we solve the in-sample problem of Section 5.1.2 to obtain $\mathbf{x}_m(\epsilon_W) \in \{0,1\}^n$, where it denotes the attacker's selected feature set under radius $\epsilon_W$.

We then assess its out-of-sample performance via 5-fold cross-validation over budgets $\mathcal{B} = \{1, \ldots, 14\}$:

$$\bar{a}_m(\epsilon_W) = \frac{1}{5\,|\mathcal{B}|} \sum_{k=1}^{5} \sum_{B \in \mathcal{B}} \text{Accuracy}\big(\mathbf{x}_m(\epsilon_W), B; k\big), \tag{30}$$

where in each fold $k$ we train on 240 perturbed datasets and test on the remaining 60. Finally, we select

$$\epsilon^*_{W,m} = \arg\min_{\epsilon_W} \; \bar{a}_m(\epsilon_W), \tag{31}$$

i.e., the radius that minimizes mean test accuracy (maximizing attack effectiveness) for each $m$.

Table 3 reports the defender's average out-of-sample test accuracy (over 5-fold cross-validation and all budgets) for the distributionally robust (DRO) and risk-receptive (DRR) attackers for each candidate $\epsilon_W$. The radius that minimizes defender accuracy for the DRO attacker is $\epsilon^*_{W,\text{RA}} = 0.1$, while for the DRR attacker $\epsilon^*_{W,\text{RR}} = 1.0$ or $1.5$. From an attacker's perspective, these results inform the selection of the Wasserstein radius ($\epsilon_W$) according to their risk appetite. For a risk-receptive attacker, employing the larger $\epsilon_W$ appears advantageous as it provides greater flexibility to identify highly impactful attacking solutions. Conversely, a risk-averse attacker may find that an overly large $\epsilon_W$ can be inefficient, as such large radius can compel an overly conservative strategy to counter an extensive range of potential distributions, paradoxically resulting in a modest increase in the defender's mean test accuracy.

| $\epsilon_W$ | DRO $k$-SIP | DRR $k$-SIP |
|---|---|---|
| 0.1 | 0.9590 | 0.9562 |
| 0.3 | 0.9597 | 0.9550 |
| 0.5 | 0.9601 | 0.9551 |
| 0.7 | 0.9599 | 0.9548 |
| 1.0 | 0.9611 | 0.9537 |
| 1.5 | 0.9612 | 0.9537 |

**Table 3**: Defender's average out-of-sample test accuracy as a function of $\epsilon_W$

## 5.2 Weighted Coverage Interdiction Problem

### 5.2.1 Instance Generation

For $k = 1$, we use instances provided by [37] for problem presented in subsection 3.2.2, where customer locations are randomly generated from Uniform $[1, 10] \times [1, 10]$ and candidate locations for the sensor is same as the generated location of the customers.

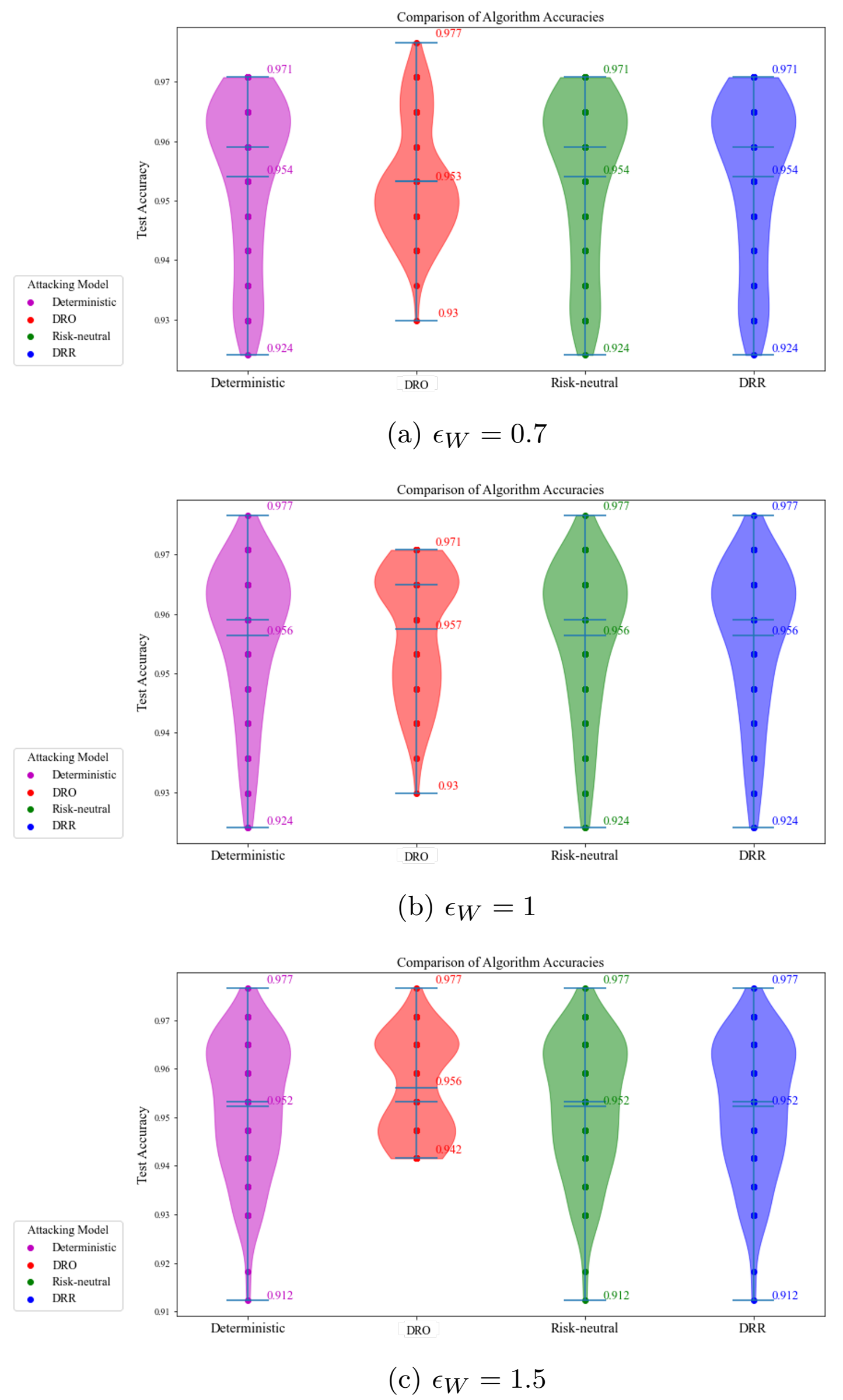

(a) $\epsilon_W = 0.7$

(b) $\epsilon_W = 1$

(c) $\epsilon_W = 1.5$

**Fig. 4**: Performance of Support Vector Classifiers (SVCs) trained on defender selected features across different problem scenarios with varying $\epsilon_W$ and fixed attacker and defender budget of 6

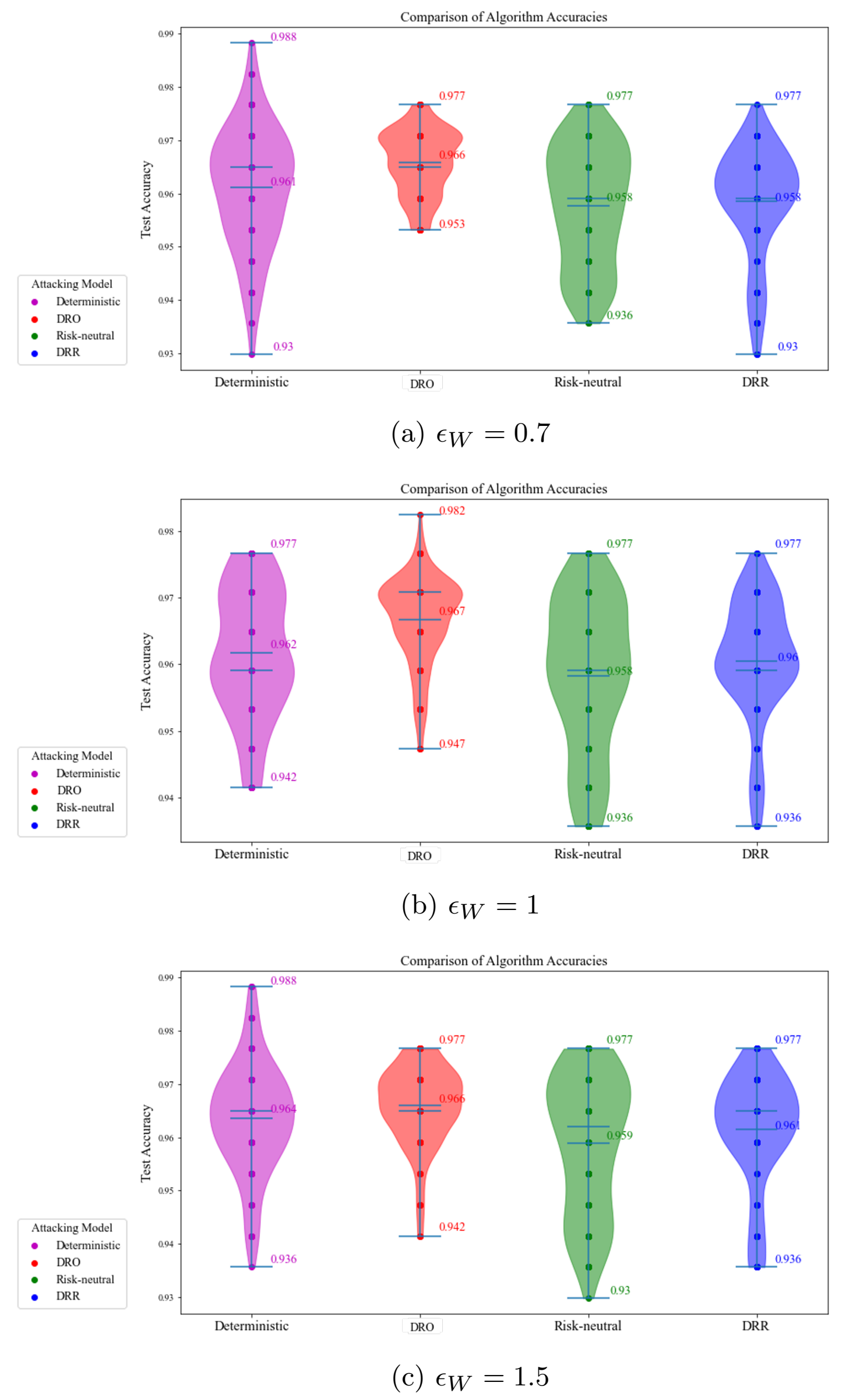

(a) $\epsilon_W = 0.7$

(b) $\epsilon_W = 1$

(c) $\epsilon_W = 1.5$

**Fig. 5**: Performance of Support Vector Classifiers (SVCs) trained on defender selected features across different problem scenarios with varying $\epsilon_W$ and fixed attacker and defender budget of 7

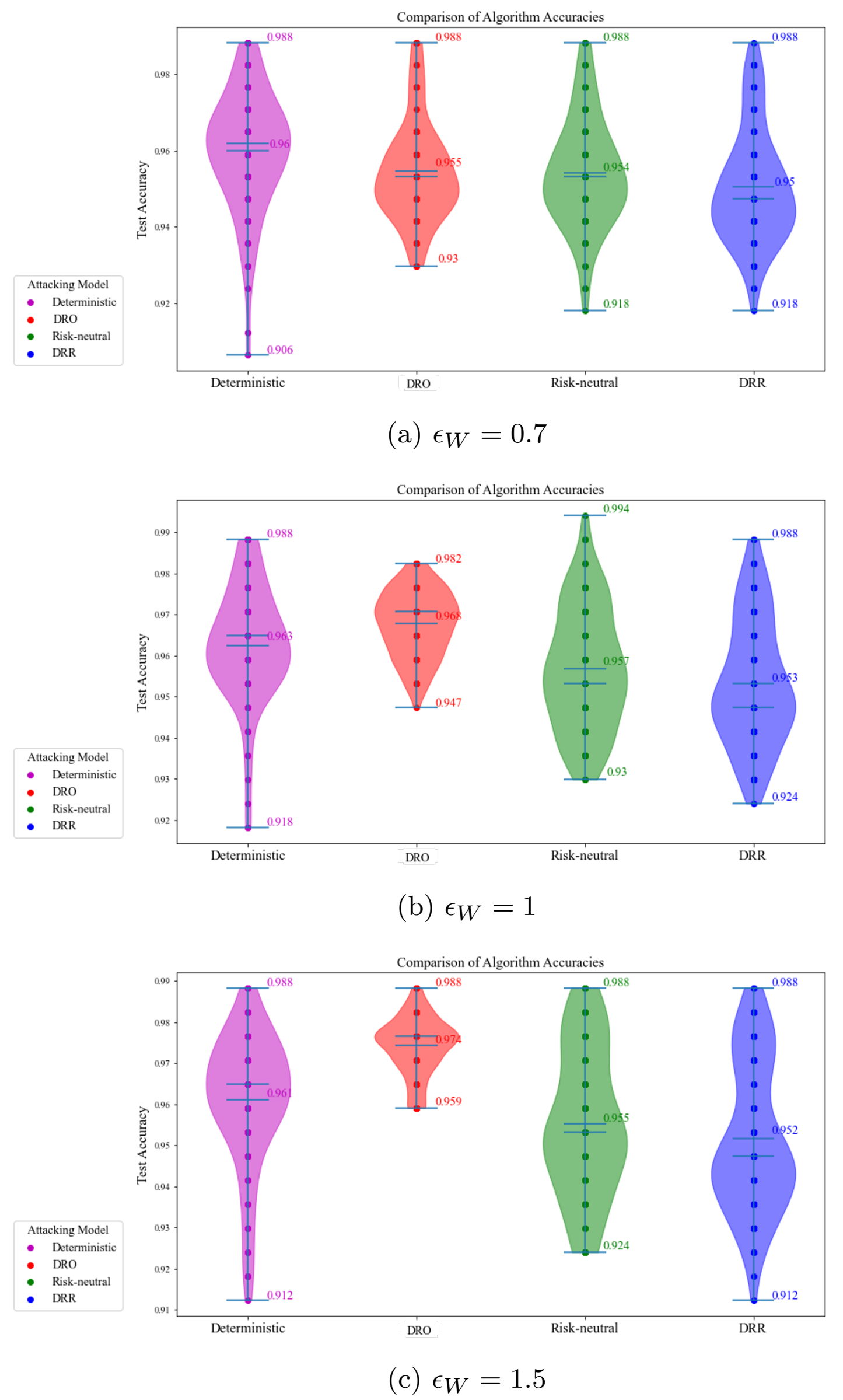

(a) $\epsilon_W = 0.7$

(b) $\epsilon_W = 1$

(c) $\epsilon_W = 1.5$

**Fig. 6**: Performance of Support Vector Classifiers (SVCs) trained on defender selected features across different problem scenarios with varying $\epsilon_W$ and fixed attacker and defender budget of 14

| $n$ | $\|\Omega\|$ | Deterministic $k$-SIP | | DRO $k$-SIP | | Risk-Neutral $k$-SIP | | DRR $k$-SIP | |
|---|---|---|---|---|---|---|---|---|---|
| | | Reward $\Phi_D(\mathbf{x}_{DT})$ | Time $(s)$ | Reward $\Phi_A(\mathbf{x}_{RA})$ | Time $(s)$ | Reward $\Phi_N(\mathbf{x}_{RN})$ | Time $(s)$ | Reward $\Phi_R(\mathbf{x}_{RR})$ | Time $(s)$ |
| 50 | 100 | 1687 | 0.37 | 1770 | 21 | 1748 | 10 | 1722 | 56 |
| 60 | | 2093 | 0.57 | 2215 | 67 | 2185 | 38 | 2147 | 252 |
| 70 | | 2597 | 1.2 | 2,717 | 361 | 2690 | 204 | 2655 | 1370 |
| 80 | | 3044 | 2.8 | 3194 | 2378 | 3164 | 993 | 3123 | 4610 |
| 90 | | 3599 | 11 | 3773 | 5722 | 3733 | 2602 | 3703 | 14690 |
| 100 | | 2925 | 37 | 3133 | 6933 | 3100 | 3906 | 3034 | 16815 |

**Table 4**: Weighted Coverage Interdiction Problem: Computational results for solving Deterministic $k$-SIP, DRO $k$-SIP, Risk-Neutral $k$-SIP and DRR $k$-SIP where $k = 1$.

For customer $j \in N$, we randomly generate associated reward $p_j$ from Uniform $[1, 100]$ and the customer $j$ is covered by sensor placed at location $i \in N$ if euclidean distance $d_{ij} \leq r$ where $r$ denotes the radius of the sensor. Instances were created for radii $r = \{1, 2\}$, maintaining uniform sensor radius within each instance. For $k = 2$, indicating the presence of two sensor types with distinct radii ($r = 1$ or $r = 2$), only one sensor of either type can be installed at each location. In our instance-generation procedure for computational studies, larger-radius sensors ($r = 2$) yield only half the reward per covered customer compared to smaller-radius sensors ($r = 1$), reflecting our designed trade-off between range and payoff. For both categories, i.e., $k = 1$ and $k = 2$, we generate a realization $\xi^\omega \in \{0, 1\}^n$, $\omega \in \Omega$, following a Bernoulli distribution with a success probability of 0.75.

### 5.2.2 Computational Results

Table 4 details the results from solving Deterministic $k$-SIP and applying our decomposition methods to solve DRO $k$-SIP, Risk-Neutral $k$-SIP and DRR $k$-SIP under the moment matching ambiguity set ($\epsilon_M = 0.05$). The "Reward" and "Time (s)" columns list the objective values and total CPU time to reach optimality (i.e., $\theta^{ub} - \theta^{lb} \leq \epsilon$), respectively. We report average over 18 instances for each row in Table 4, except for $n = 100$, where the results are averaged over 8 instances. For each $n$, we assume that budget for the attacker and defender is $0.1 \times n$. As expected, defender's covered rewards for DRR $k$-SIP are consistently less than those for Risk-Neutral $k$-SIP and DRO $k$-SIP. This suggests that DRR $k$-SIP could serve as a strategic tool for attackers who are willing to take a risk, providing them with a quantitative measure of how much they can potentially diminish the defender's objective of covered demand. From the vulnerability analysis standpoint for the defender, the DRR $k$-SIP identifies the most severe attacking strategies, enabling the identification of critical locations, which can significantly reduce captured reward if compromised. This solution allows the defender to plan fortification of locations against potential attacks obtained from DRR $k$-SIP. This approach is reasonable, as defenders cannot anticipate the risk preferences of attackers, requiring them to prepare for a risk-receptive attackers as well. Conversely, for attackers who prefer to avoid a risk, the DRO $k$-SIP offers a conservative estimate of the maximum reward the defender could achieve, thereby informing a risk-averse strategy. Note that defender's $\Phi_D(\mathbf{x}_{DT})$ are less than those of the other three problems, reflecting the naive assumption in the deterministic problem that all attacks are successful. This assumption, although optimistic from the attacker's perspective, does not realistically capture the uncertain nature of real-world scenarios, potentially misleading the attacker in determining the best attacking solution.

| $n$ | $\Phi_N(\mathbf{x}_{DT})$ | $\Phi_N(\mathbf{x}_{RN})$ | VSS | # of instances VSS>0 | $\Phi_A(\mathbf{x}_{DT})$ | $\Phi_A(\mathbf{x}_{RA})$ | VAS | # of instances VAS>0 | $\Phi_R(\mathbf{x}_{DT})$ | $\Phi_R(\mathbf{x}_{RR})$ | VRS | # of instances VRS>0 |
|---|---|---|---|---|---|---|---|---|---|---|---|---|
| 50 | 1760.7 | 1748.0 | 25.3 | 9 | 1787.6 | 1770.5 | 25.7 | 12 | 1732.2 | 1722.7 | 21.5 | 8 |
| 60 | 2197.5 | 2185.5 | 19.6 | 11 | 2234.0 | 2215.4 | 24.0 | 14 | 2158.0 | 2147.9 | 20.1 | 9 |
| 70 | 2707.1 | 2690.2 | 20.3 | 15 | 2740.1 | 2717.2 | 24.2 | 17 | 2671.8 | 2655.6 | 19.5 | 15 |
| 80 | 3182 | 3164.5 | 23.2 | 14 | 3220.0 | 3194.6 | 35.1 | 13 | 3138.8 | 3123.7 | 20.8 | 13 |
| 90 | 3770.5 | 3743.1 | 30.9 | 16 | 3806.5 | 3773.2 | 49.9 | 12 | 3723.5 | 3703.8 | 32.4 | 11 |
| 100 | 3156.8 | 3100.7 | 64.0 | 7 | 3330.5 | 3133.4 | 197.5 | 8 | 3038.1 | 3034.0 | 32.5 | 1 |

**Table 5**: Average values of VSS, VAS and VRS

To underscore the importance of accounting for uncertainty in the optimization models, we measure Value of Stochastic Solution (VSS) that is computed as follows.

$$\text{VSS} = \Phi_N(\mathbf{x}_{DT}) - \Phi_N(\mathbf{x}_{RN}),$$

where $\Phi_N(\mathbf{x}_{DT})$ returns the defender's expected covered reward when attacker's adheres to $\mathbf{x}_{DT}$ (deterministic optimal solution) and there is no distributional ambiguity. Additionally, we introduce Value of Distributionally Risk Averse Solution (VAS) and Value of Distributionally Risk Receptive Solution (VRS) as follows.

$$\text{VAS} = \Phi_A(\mathbf{x}_{DT}) - \Phi_A(\mathbf{x}_{RA}) \text{ and } \text{VRS} = \Phi_R(\mathbf{x}_{DT}) - \Phi_R(\mathbf{x}_{RR}),$$

By definition, VSS, VAS and VRS are non-negative as $\mathbf{x}_{RN}, \mathbf{x}_{RA}$ and $\mathbf{x}_{RR}$ are all optimal solutions and $\mathbf{x}_{DT}$ is a feasible solution. In Table 5, we report average of VSS only for instances where VSS $> 0$, out of 18 instances for each $n$ , except for $n = 100$, and similarly, average VAS and VRS are reported for instances where VAS $> 0$ and VRS $> 0$.

From Table 5, when $n = 100$, VAS is 197.5 and this suggests that if the attacker naively relies on the deterministic optimal solution $\mathbf{x}_{DT}$ despite the presence of uncertainty, they miss the opportunity to diminish the defender's reward by 197.5 units, which could have been achieved by solving DRO $k$-SIP and opting $\mathbf{x}_{RA}$. This interpretation holds when VSS or VRS is greater than zero.

As shown in Figure 7, the values of VAS, VSS, and VRS at $n = 50$ and $n = 60$ are similar, with a modest dip at $n = 60$ compared to $n = 50$. From $n = 70$ onward, all three measures increase steadily with $n$, and VAS remains above both VSS and VRS for every subsequent $n$. These observations suggest that, as the number of candidate locations and budgets grow, stochastic solutions $\mathbf{x}_{RA}$, $\mathbf{x}_{RN}$, and $\mathbf{x}_{RR}$ tend to offer gradually larger advantages over the deterministic attacker's solution.

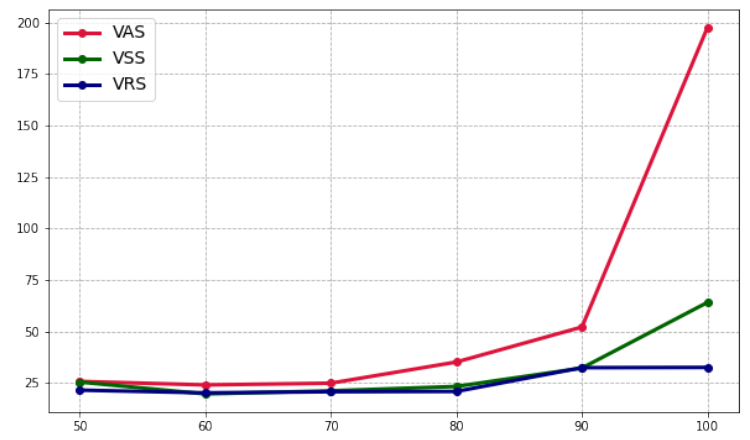

**Fig. 7**: Comparison of VSS, VAS and VRS

*Weighted Coverage Interdiction Problem with* $k = 2$. Table 6 details the results for bi-submodular function. The results are averaged over 10 instances for each row in Table 4, and for each $n$, we assume that budget for the attacker and defender is $0.1 \times n$. Additionally, we set the number of scenarios $|\Omega| = 30$. We observe that defender's optimal captured rewards for DRR 2-SIP is less than those of Risk-neutral 2-SIP and DRO 2-SIP. Again, this results aligns with the fact that DRR $k$-SIP provides the solution where the attacker can degrade the defender's objective to the greatest extent, making it useful for vulnerability analysis from the defender's perspective. Conversely, the defender's capture rewards for DRO 2-SIP are greater than those of Risk-neutral 2-SIP and DRR 2-SIP, as it provides the conservative strategy towards the distributional ambiguity for the attacker.

| $n$ | $\|\Omega\|$ | Deterministic $k$-SIP | | DRO $k$-SIP | | Risk-Neutral $k$-SIP | | DRR $k$-SIP | |
|---|---|---|---|---|---|---|---|---|---|
| | | Reward $\Phi_D(\mathbf{x}_{DT})$ | Time $(s)$ | Reward $\Phi_A(\mathbf{x}_{RA})$ | Time $(s)$ | Reward $\Phi_N(\mathbf{x}_{RN})$ | Time $(s)$ | Reward $\Phi_R(\mathbf{x}_{RR})$ | Time $(s)$ |
| 10 | | 83.0 | 0.09 | 84.6 | 0.1 | 84.4 | 0.1 | 84.2 | 0.2 |
| 20 | | 310 | 0.2 | 323 | 1.3 | 321 | 1.4 | 318 | 1.5 |
| 30 | | 730 | 0.7 | 842 | 8.8 | 823 | 8.6 | 802 | 14 |
| 40 | 30 | 933 | 6.2 | 1045 | 76 | 1030 | 85 | 1014 | 129 |
| 50 | | 1278 | 9.5 | 1452 | 142 | 1428 | 150 | 1408 | 238 |
| 60 | | 2013 | 2550 | 2160 | 7053 | 2142 | 8884 | 2124 | 16870 |
| 70 | | 2732 | 1425 | 2963 | 4341 | 2932 | 5033 | 2898 | 10490 |

**Table 6**: Weighted Coverage Interdiction Problem: Computational results for solving Deterministic k-SIP, DRO k-SIP, Risk-Neutral k-SIP and DRR k-SIP where $k = 2$.

# 6 Conclusion

To address submodular optimization in adversarial and uncertain environment, we introduced Distributionally Robust $k$-Submodular Interdiction Problem (DRO $k$-SIP) and Distributionally Risk-Receptive $k$-Submodular Interdiction Problem (DRR $k$-SIP) and presented exact solution approaches for them. We conducted computational experiments on instances of Feature Selection Interdiction Problem (FSIP) and Multi-type Sensor Coverage Interdiction Problem. Note that feature selection is a key concept in machine learning, underscoring the importance of these results for practical applications. We analyzed the solutions obtained from both problems from each of the decision maker's perspective with varying levels of risk-appetite. In FSIP, the optimal solution from DRO 1-SIP demonstrated robust feature removal strategy which is effective from the attacker's risk-averse (or conservative) perspective. Conversely, solution from DRR 1-SIP identified critical features whose removal reduce the quality of data the most from the attacker's risk-receptive (or defender's risk-averse) perspective. In general, the DRO $k$-SIP seeks to determine the optimal expected value of the defender's objective function under the worst probability distribution from the attacker's perspective, providing a robust strategy suitable for risk-averse attackers. In contrast, the DRR $k$-SIP offers insights into effective strategies for risk-taking attackers, thereby identifying critical vulnerabilities in the defender's system.

**Author contributions** The authors contributed equally to the manuscript.

**Funding** This research is partially funded by National Science Foundation Grant CMMI-1824897 and 2034503, and Commonwealth Cyber Initiative grants which are gratefully acknowledged.

**Data availability** The instances used for computational studies in this paper will be made available in "Submodular-Interdiction-Game" folder at https://github.com/Bansal-ORGroup/.

# Declarations

- **Conflict of interest** The authors report that there are no Conflict of interest to declare.